\documentclass{article}
\usepackage{amsmath}
\usepackage{amsthm}
\usepackage{amssymb}
\usepackage{amsfonts}

\usepackage{enumitem}
\usepackage{tikz}
\usepackage{pgfplots}
\usetikzlibrary{arrows.meta}

\usepackage[heads=littlevee]{diagrams}

\usepackage[labelsep=newline]{caption}

\numberwithin{equation}{section}
\numberwithin{figure}{section}

\newtheorem{lem}[equation]{Lemma}
\newtheorem{thm}[equation]{Theorem}

\newcommand{\labelthis}[1]{\addtocounter{equation}{1}\tag{\theequation}\label{#1}}

\renewcommand*{\bar}{\overline}
\newcommand{\bra}[1]{\left(#1\right)}
\newcommand{\abs}[1]{\left|#1\right|}
\newcommand{\norm}[1]{{\left\Vert#1\right\Vert}}
\newcommand{\Set}[2]{\left\{\; #1 \;\middle|\; #2 \;\right\}}
\newcommand{\vt}[2]{\begin{pmatrix} #1 \\ #2 \end{pmatrix}}
\newcommand{\Partial}[2]{\frac{\partial #1}{\partial #2}}
\newcommand{\cross}[4]{\left[#1,#2;\, #3, #4\right]}
\newcommand{\cji}[1]{\overline{#1}^{-1}}
\newcommand{\diffeo}{\cong}

\DeclareMathOperator{\atan}{atan}

\DeclareMathOperator{\asinh}{asinh}

\DeclareMathOperator{\Wind}{Wind}

\DeclareMathOperator{\Mat}{Mat}

\DeclareMathOperator*{\tr}{tr}

\DeclareMathOperator{\pp}{pp.}

\DeclareMathOperator{\Real}{Re}
\DeclareMathOperator{\Imag}{Im}

\newcommand{\iu}{i}

\newcommand{\N}{\mathbb{N}}
\newcommand{\Z}{\mathbb{Z}}
\newcommand{\Q}{\mathbb{Q}}
\newcommand{\R}{\mathbb{R}}

\newcommand{\C}{\mathbb{C}}
\newcommand{\RP}{\mathbb{R}\text{P}}
\newcommand{\CP}{\mathbb{C}\text{P}}

\renewcommand{\S}{\mathbb{S}}

\newcommand{\SU}{\mathrm{SU}}
\newcommand{\SL}{\mathrm{SL}}
\newcommand{\SO}{\mathrm{SO}}
\newcommand{\su}{\mathfrak{su}}

\begin{document}
% FOR ARXIV
\title{The space of equivariant harmonic tori in the $3$-sphere}
\author{Emma Carberry and Ross Ogilvie}
\date{}
\maketitle

\begin{abstract}
In this paper we give an explicit parametrisation of the moduli space of equivariant harmonic maps from a $2$-torus to $\S^3$. 
As Hitchin proved, a harmonic map of a $2$-torus is described by its spectral data, which consists of a hyperelliptic curve together with a pair of differentials and a line bundle. The space of spectral data is naturally a fibre bundle over the space of spectral curves.
For homogeneous tori the space of spectral curves is a disc and the bundle is trivial.
For tori with a one-dimensional invariance group, we enumerate the path connected components of the space of spectral curves and show that they are either `helicoids' or annuli, and that they densely foliate the parameter space. The bundle structure of the moduli space of spectral data over the annuli components is nontrivial.
In the two cases, the spectral data require only elementary and elliptic functions respectively and we give explicit formulae at every stage.
Homogeneous tori and the Gauss maps of Delaunay cylinders are used as illustrative examples.
\end{abstract}

\section{Introduction}\label{sec:Introduction}

In this paper we consider harmonic maps from a $2$-torus to the unit $3$-sphere $\S^3$ that are equivariant with respect to the action of $\SO(4)$.
Hitchin~\cite{Hitchin1990} gave a correspondence between these harmonic maps and hyperelliptic curves with additional data. A hyperelliptic curve that occurs in this correspondence is called a spectral curve. 
The genus of the spectral curve is a natural invariant of the harmonic map; homogeneous tori are exactly those whose spectral curve is genus zero, whereas those with only a one-dimensional invariance group have spectral genus one.
The additional data of these spectral curves include marked points, differentials of the second kind, and a line bundle. 
Given a spectral curve one can choose the differentials from a $\mathbb{Z}^2$ lattice and the line bundle from a real subvariety of the Jacobian. Therefore the main difficulty is to describe the moduli space of spectral curves.

It was already known that the space $\mathcal{S}_0$ of spectral curves of genus zero is an open disc~\cite[Section 9]{Hitchin1990}.
For the space $\mathcal{S}_1$ of spectral curves of genus one, we are able to enumerate the path connected components and show that they are two-dimensional surfaces, either contractible or annuli. Moreover, by constructing an appropriate coordinate system we foliate the parameter space of hyperelliptic curves and exhibit $\mathcal{S}_1$ as a dense collection of leaves.
The bundle structure of the space $\mathcal{M}_1$ of spectral data over the contractible components is obviously trivial, but over the annuli components it is nontrivial. We will explain its structure in terms of a monodromy action on the $\Z^2$ lattice of differentials.

Our investigation was motivated by~\cite{Kilian2015}, which examines constant mean curvature (CMC) tori in $\S^3$. 
Kilian, Schmidt, and Schmitt establish that generically the spectral curves of CMC tori come in one-dimensional families. The theories of harmonic and CMC surfaces overlap in the sense that a conformally parametrised surface is harmonic if and only if it has mean curvature zero.
Another point of connection, due to Ruh-Vilms~\cite{Ruh1970}, is that the Gauss maps of conformally parametrised CMC surfaces in $\R^3$ are non-conformal harmonic maps to a great-sphere in $\S^3$.
In contrast to CMC tori, for harmonic tori it is known that the moduli space of spectral curves of any genus is generically two-dimensional~\cite{Carberry2019} (and in particular the spaces of spectral curves with genus zero, one, and two are smooth surfaces). Thus there is a broader range of possibilities for the topology of the space and how it is embedded within the space of hyperelliptic curves, and different methods were required.

A good illustration of the space $\mathcal{S}_1$ of spectral curves of genus one is given by its cross sections. In our context, there is a natural parameter space in which $\mathcal{S}_1$ lives, but it is not simply connected. We will construct coordinate charts on the universal cover of the parameter space of the form
\[
\left\{ \bra{p,q,k,\tilde{X}} \in I \times \R \times (0,1) \times \R \right\},
\labelthis{eqn:parameter demo}
\]
where $I\subset\R^+$ is an interval. The preimage of $\mathcal{S}_1$ is given in these coordinates by $p\in\Q^+$ and $q\in\Q$. 
This decomposes the preimage into path connected components. 
We must then push this decomposition back down into the space of curves.
The group of covering transformations is generated by $\tilde{\lambda}$ which acts as
\[
\tilde{\lambda} : \bra{p,q,k,\tilde{X}} \mapsto \bra{p,q + (p-1),k,\tilde{X} + \pi}.
\]
% Thus $\tilde{\lambda}$ maps $\mathcal{\tilde{A}}(p,q)$ to $\mathcal{\tilde{A}}(p,q + p-1)$. 
% If $p$ and $q$ are rational, so too is $q + (p-1)$ and hence the group $\mathcal{G} = \Z\langle \tilde{\lambda}\rangle$ of covering transformations restricts to give a group action on the preimage of the space of spectral curves.
Since both $p\in\Q^+$ and $k$ are unchanged by the group action, let us take a cross section $(q,\tilde{X}) \in \R^2$, shown in Figure~\ref{fig:level set quotient}. When we quotient the universal cover by $\langle\tilde{\lambda}\rangle$, for $p\neq1$ the $(q,\tilde{X})-$plane becomes a cylinder foliated by helices $q=\text{const}$. If we to allow $k$ to vary, imagine the cylinder instead as a solid cylinder with the centre axis removed. The helices extend to `helicoids' in this coordinate system (but we do not use a metric on the parameter space, the apparent constant pitch of this surface is not meaningful, and we use quotation marks to remind ourselves of this fact). The `helicoids' where $q$ is rational correspond to spectral curves. These are plotted in the natural coordinates of the parameter space in Figure~\ref{fig:p05 plot}.

\begin{figure}
    % \includestandalone[width=\textwidth]{tikz/universal_quotient}
\resizebox{\textwidth}{!}{
\begin{tikzpicture}
    \clip (-3.5,-3.5) rectangle (11.5,3.5);
    
    % \draw[color=gray,step=1.0,dotted] (-3.3,-3.3) grid (3.3,3.3);
    
    \draw[->] (-3.3,-3.3)--(3.3,3.3) node[below=3pt]{$q$};
    \draw[->] (3.3,-3.3)--(-3.3,3.3) node[below=3pt]{$\tilde{X}$};
    \draw[->,color=blue] (0,0)--(0,2) node[right]{$\tilde{\lambda}$};
    
    % \draw[->] (0,0.7) arc (90:135:0.7) node[above=3pt]{$\theta$};
    
    \foreach \x in {1,...,12}
        \draw[color=gray, thin] (-3.3, 3.3-0.5*\x) -- (3.3-0.5*\x, -3.3);
    \foreach \x in {1,...,12}
        \draw[color=gray, thin] (-3.3+0.5*\x, 3.3) -- (3.3, -3.3+0.5*\x);
    
    \draw (5,0) ellipse (0.5 and 1);
    \draw[->,color=blue] (7.85,0)--(7.85,0.01) node[right]{$\tilde{\lambda}$};
    \draw[color=blue] (7.5,1) arc (40:-40:1.56);
    \draw (10,1) arc (60:-60:1.1547);
    \draw (5,1) -- (10,1);
    \draw (5,-1) -- (10,-1);

    \draw[color=gray, thin] (6.9,-1) arc (15:67:2.75);
    \draw[color=gray, thin] (6.5,-1) arc (15:55:2.75);
    \draw[color=gray, thin] (6.1,-1) arc (15:42:2.75);
    \draw[color=gray, thin] (5.7,-1) arc (15:29:2.75);
    \foreach \x in {0,...,6}
    \draw[color=gray, thin] (5.2+0.4*\x,1) arc (80:15:2.75);
    \draw[color=gray, thin] (8,1) arc (80:17:2.75);
    \draw[color=gray, thin] (8.4,1) arc (80:24.5:2.75);
    \draw[color=gray, thin] (8.8,1) arc (80:35:2.75);
    \draw[color=gray, thin] (9.2,1) arc (80:50:2.75);
    \draw[color=gray, thin] (9.6,1) arc (80:65:2.75);
    
    % \fill (2.6,2.3) circle (0.05) node[left,color=blue]{$0$} node[right,color=black]{$z_0$};
    % \fill (-2.6,2.4) circle (0.05) node[left, color=blue]{$\infty$} node[right,color=black]{$-\bar{z}_0$};
    %
    % \fill (0,2) circle (0.05) node[left]{$f(1)$};
    % \fill (0,-2.4) circle (0.05) node[left]{$f(-1)$};
    %
    % \fill (0,0) circle (0.05) node[above=5pt, right, color=blue]{$\mu$} node[below=8.5pt, left,color=black]{$0$};
    % \fill (1,0) circle (0.05) node[above, color=blue]{$\alpha$} node[below=1.7pt,color=black]{$1$};
    % \fill (2.1,0) circle (0.05) node[above, color=blue]{$\beta$} node[below,color=black]{$k^{-1}$};
    % \fill (-1,0) circle (0.05) node[above, color=blue]{$\bar{\alpha}^{-1}$} node[below=1.7pt,color=black]{$-1$};
    % \fill (-2.1,0) circle (0.05) node[above, color=blue]{$\bar{\beta}^{-1}$} node[below,color=black]{$-k^{-1}$};
    
\end{tikzpicture}
}
    \caption{On the left, a cross-section of the universal cover of the parameter space where $p$ and $k$ have been fixed. The translation $\tilde{\lambda}$ is in blue. The grey lines are the lines of constant $q$.\\
    On the right, the result of taking the quotient by $\langle \tilde{\lambda}\rangle$. The plane has been rolled into a cylinder, and the level sets are now helices.\label{fig:level set quotient}}
\end{figure}
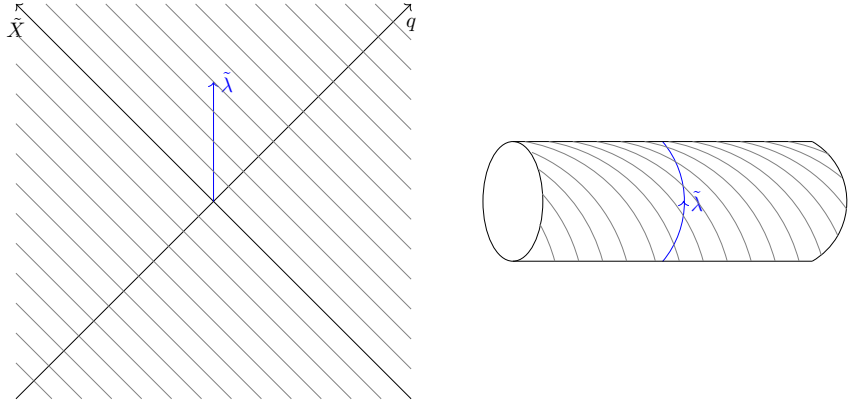
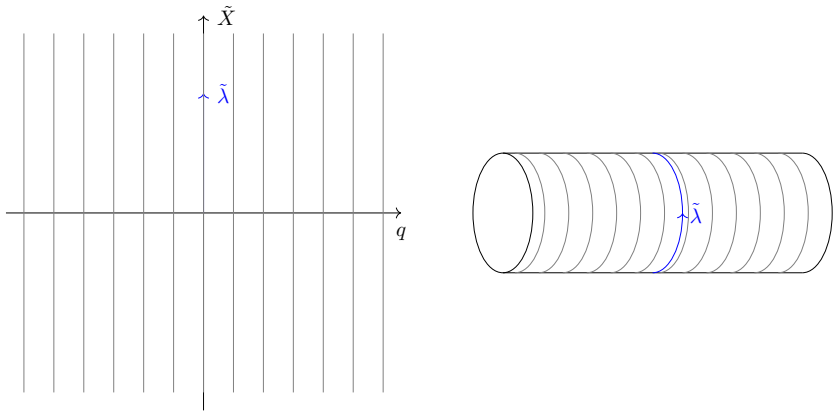
\begin{figure}
    % \includestandalone[width=\textwidth]{tikz/universal_quotient}
\resizebox{\textwidth}{!}{
\begin{tikzpicture}
    \clip (-3.5,-3.5) rectangle (11.5,3.5);
    
    % PLANE ON LEFT
    \draw[->] (-3.3,0)--(3.3,0) node[below=3pt]{$q$};
    \draw[->] (0,-3.3)--(0,3.3) node[right=3pt]{$\tilde{X}$};
    \draw[->,color=blue] (0,0)--(0,2) node[right=3pt]{$\tilde{\lambda}$};
    
    \foreach \x in {0,...,12}
        \draw[color=gray, thin] (-3+0.5*\x,-3) -- (-3+0.5*\x,3);
    
    % CYLINDER ON RIGHT
    \draw (5,0) ellipse (0.5 and 1);
    \draw (5,1) -- (10,1);
    \draw (5,-1) -- (10,-1);
    \draw (10,-1) arc (-90:90:0.5 and 1);

    \draw[->,color=blue] (7.5,-1) arc (-90:0:0.5 and 1) node[right]{$\tilde{\lambda}$};
    \draw[color=blue] (8,0) arc (0:90:0.5 and 1);
    
    \foreach \x in {0,...,11}
        \draw[color=gray, thin] (5.2+0.4*\x,-1) arc (-90:90:0.5 and 1);
    \end{tikzpicture}
}
    
    \caption{When $p=1$, the translation $\tilde{\lambda}$ fixes $q$. When taking the quotient, these level sets become circles.\label{fig:level set quotient p1}}
\end{figure}

One can understand the role of $p$ as determining the angle between the translation $\tilde{\lambda}$ and the coordinate axes. As $p$ tends to $1$, the translation aligns more closely with the $q$-axis. When $p=1$ lines of constant $q$ are in the direction of translation and close up into circles under the group action. We think of this exceptional case $p=1$ as the transition from left-handed to right-handed `helicoids'. The cross-section diagram is shown in Figure~\ref{fig:level set quotient p1} or in full in Figure~\ref{fig:p1 plot}.

A precise statement of our results is the following. In {\bf Theorem~\ref{thm:topology_curves}} we prove that the subset of the parameter space with $p\neq 1$ constant is diffeomorphic to 
\[
\left\{ \bra{[q],k,\tilde{X}} \in \R/(p-1)\Z \times (0,1) \times \R \right\}.
\]
The points where $[q] \in \Q/(p-1)\Z$ correspond to spectral curves. This is a union of `helicoids'.
For $p=1$, by contrast, {\bf Theorem~\ref{thm:topology_curves_p1}} shows that the subset is 
\[
\left\{ \bra{q,k,\left[\tilde{X}\right]} \in \R \times (0,1) \times \R/\pi\Z \right\}.
\]
This time the spectral curves are the subset where $q\in\Q$ and is a union of annuli.

Although there is a classification of harmonic and CMC maps given by~\cite{Hitchin1990,Pinkall1989,Bobenko1991} in terms of algebraic curves, many questions remain. 
On one hand, we can study local deformations of a map that preserve properties such as being harmonic or CMC. The idea of spectral curves arose in the study of classical PDE systems, such as KdV, and the spectrum of their operators~\cite{Flaschka1980,Krichever1995}. Linearisations of the differential equation are understood in the correspondence as changing the line bundle without changing the spectral curve, and thus they are also known as isospectral deformations. More recent work, for example~\cite{Kilian2010} and~\cite{Carberry2019}, use Whitham deformations, which allow for the spectral curve and differentials to vary as well. Such results often establish the existence of smooth components in the moduli space and compute their dimensions.
On the other hand, there are results that show spectral curves to be dense in some parameter space, for example~\cite{Carberry2016a,Carberry2016}.
Our results here address both questions simultaneously, showing explicitly that the space of harmonic maps is two-dimensional but also that it densely foliates a space of hyperelliptic curves.
More details about the spectral curve correspondence may be found in the expository papers~\cite{McIntosh2008} and~\cite{Carberry2013a}.

The results of this paper indicate at least two interesting avenues for further investigation. Most obviously, one could attempt this programme for higher genus spectral curves using theta functions and a braid space, though one would probably have to limit oneself to non-singular spectral curves. This is due to what we see in Figures~\ref{fig:p05 plot}--\ref{fig:p1 plot}, that the space of spectral curves has peculiar behaviour near a singular curve; it spirals and the different path components converge to a single point. It is unclear whether this is particular to this low genus case or typical of the space of spectral curves generally.

On a different part of the boundary, in a suitable limit spectral curves accumulate a different kind of singularity and can be normalised to lower genus curve. This has been carried out for CMC tori in $\S^3$ with spectral genus one in~\cite{Kilian2015}, and in general for spectral data with only one differential in~\cite{Hauswirth2017arxiv}. The formulae in this paper allow one to observe the differentials during this procedure and to begin to connect together the spaces of spectral curves with different genera. The preliminaries of this process is covered in~\cite[Section~4.2]{Ogilvie2017}. 

We now give a section-by-section outline of the paper. Section~\ref{sec:Genus Zero} deals with harmonic maps whose spectral curve have genus zero. This case was already well-understood, but we use it as an introduction to spectral data. Starting with explicit formulae of the harmonic maps, we construct their spectral data using a family of flat connections~\eqref{eqn:flat connection translation} and their holonomy.
We also discuss how the geometry of a harmonic map is altered as one varies its spectral data.

The main results of the paper are contained in Section~\ref{sec:Genus One}. In this section we begin with a parameter space for genus one hyperelliptic curves and then successively impose conditions that a curve must satisfy to be a spectral curve. This gives explicit forms for the spectral data in terms of inverse elliptic functions.
The section culminates in Theorems~\ref{thm:topology_curves} and~\ref{thm:topology_curves_p1} which describe the topology of the space of genus one spectral curves, enumerates the the path connected components, and shows that they are dense in the space of hyperelliptic curves. It concludes with an explanation of the bundle structure of the space $\mathcal{M}_1$ of spectral data over this space $\mathcal{S}_1$ of spectral curves.

The final section takes a symmetry used in Lemma~\ref{lem:T_graph} and generalises it to spectral curves of any genus. We compare how this symmetry acts for genus zero and for genus one. In the latter case, we connect this to the Gauss maps of Delaunay cylinders.

The authors would like to thank Prof John Rice for his feedback as the material was developed. The second author also wishes to thank Dr Prof Martin Schmidt for discussions about spectral genus one CMC tori.

We close the introduction by recalling a theorem of Hitchin~\cite[Theorem~8.1]{Hitchin1990}. This provides a total characterisation of spectral data of a harmonic map of a torus in $\S^3$. We therefore search for spectral curves within the space of hyperelliptic curves by examining whether a curve can satisfy all of these conditions.

\begin{thm}[Hitchin]\label{thm:Hitchin}
A tuple $(\Sigma,\Theta^1,\Theta^2,E)$ consisting of an algebraic curve $\Sigma$ with equation $\eta^2 = P(\zeta)$, a pair of differentials $\Theta^1, \Theta^2$, and a line bundle $E$ is the \emph{spectral data} of a harmonic map from the 2-torus to $\S^3$ if and only if the following conditions hold.
\begin{enumerate}[label=(P.\arabic*)]
\item\label{P:real curve} Real spectral curve: $P(\zeta)$ is a real section of $\mathcal{O}(2g+2)$ with respect to the real structure induced by $\rho : \zeta \mapsto \bar{\zeta}^{-1}$.
\item\label{P:no real zeroes} No real zeroes: $P(\zeta)$ has no zeroes on the unit circle $\S^1\subset\CP^1$.
% \item\label{P:simple zeroes} Simple zeroes: $P(\zeta)$ has only simple zeroes.
\item\label{P:poles} Poles: The differentials have double poles with no residues at $\pi^{-1}\{0,\infty\}$, where $\pi$ is the projection $\pi: \Sigma \to \CP^1$, but are otherwise holomorphic.
\item\label{P:symmetry} Symmetry: The differentials satisfy $\sigma^* \Theta = - \Theta$, for $\sigma$ the hyperelliptic involution of $\pi$.
\item\label{P:reality} Reality: The differentials satisfy $\rho^* \Theta = - \bar{\Theta}$.
\item\label{P:imaginary periods} Imaginary Periods: The differentials have purely imaginary periods.
\item\label{P:periods} Periods: The periods of the differentials lie in $2\pi\iu\Z$.
\item\label{P:closing} Closing conditions: Suppose that $\gamma^+$ is a path in $\Sigma$ connecting the two points $\pi^{-1}\{1\}$, and $\gamma^-$ connects the two points $\pi^{-1}\{-1\}$. The differentials satisfy
\[
\int_{\gamma^+}\Theta \in 2\pi\iu\Z,\;\; 
\int_{\gamma^-}\Theta \in 2\pi\iu\Z.
\]
\item\label{P:linear independence} Linear independence: The principal parts of the differentials $\Theta^1$ and $\Theta^2$ are real linearly independent.
\item\label{P:quaternionic} Quaternionic: $E^*$ is a line bundle of degree $g+1$ that is quaternionic with respect to the involution $\rho \circ \sigma$.
\end{enumerate}
\end{thm}

%%%%%%%%%%%%%%%%%%%%%%%%%%%%%%%%%%%%%%%%%%%%%%%%%%%%%%%%%%%%%%%%%%%%%%%%
%%%%%%%%%%%%%%%%%%%%%%%%%%%%%%%%%%%%%%%%%%%%%%%%%%%%%%%%%%%%%%%%%%%%%%%%
%%%%%%%% GENUS ZERO SECTION
%%%%%%%%%%%%%%%%%%%%%%%%%%%%%%%%%%%%%%%%%%%%%%%%%%%%%%%%%%%%%%%%%%%%%%%%
%%%%%%%%%%%%%%%%%%%%%%%%%%%%%%%%%%%%%%%%%%%%%%%%%%%%%%%%%%%%%%%%%%%%%%%%
\section{Spectral Curves of Genus Zero}\label{sec:Genus Zero}

This section gives a description of the space $\mathcal{M}_0$ of harmonic maps of tori whose spectral curve has genus zero: it is a discrete fibre bundle over the space of the open disc. This process serves as a guide for the description of genus one spectral curves in the subsequent section. It is an instructive case as it requires only elementary functions, and it is straightforward to give explicit formulae for all such harmonic maps $g$ and spectral data $(\Sigma,\Theta^1,\Theta^2,E)$. 
% The spectral data also serves as an example of the spectral curve construction that was outlined in Section~\ref{sec:Introduction}. 

% The derivations within this chapter may be divided into three parts. In the first part, we start with equations~\eqref{eqn:Hit1.7} and find all translation invariant solutions. Among these solutions, we determine which correspond to harmonic maps from a torus by forcing a periodicity constraint. We then write an explicit formula for each harmonic map, and bring it into a standard form by applying rotations. In the second part, we take these harmonic maps and work through the steps of Section~\ref{sec:construction} to produce the associated spectral data. Further calculations give rise to formulae for the infinitesimal deformations in terms of derivative of the branch point, and an expression of the energy of the harmonic map.

To begin, we will construct tori in $S^3 = SU(2)$. Any great circle through the identity is given by the exponential map of a line in the Lie algebra. These circles are also one-parameter subgroups. If we take one circle and translate it by another one-parameter subgroup the result is a torus in $SU(2)$. Choose two vectors $X,Y\in \su_2$. If we write $w = w_R + \iu w_I \in \C$ then the following is a map from $\C$ to $SU(2)$ whose image is a torus. The factors of $4$ have been chosen to simplify later calculations. 
\[
g(w) = \exp( -4 w_R X ) \exp( 4 w_I Y ).
\labelthis{eqn:genus zero simple map}
\]
The correspondence between spectral data and harmonic maps does not distinguish between maps that differ by an $\SO(4)$ rotation of $\S^3$. We may use this freedom to simplify $X$ and $Y$. 
% We will make both off-diagonal to simplify a computation of eigenvectors below. 
Further, by rescaling $w$ we may change the lengths of the vectors $X$ and $Y$. Without loss of generality, for some $x \in \R^+$ and $\delta\in (0,\pi)$, we may take 
\[
X = \begin{pmatrix}
0 & 1 \\ -1 & 0
\end{pmatrix}, \qquad
Y = x\begin{pmatrix}
0 & e^{\iu \delta} \\ -e^{-\iu \delta} & 0
\end{pmatrix}.
% \labelthis{eqn:def delta}
\]

Using the standard inner product on $\su_2$, which is same as the endowing $\S^3$ with the usual unit sphere metric, we may consider $\delta$ as the angle between $X$ and $Y$, and $x$ as the ratio of their lengths. The image of $g$ is determined solely by the angle $\delta$. The image is not determined uniquely by $\delta$ though; the tori for $\delta$ and $\pi-\delta$ are the same. Indeed, if one reflects $X$ to $-X$, the angle between $-X$ and $Y$ is the supplement of $\delta$. However, the two tori carry the opposite orientation and so are distinguishable.

Another way to understand these images is to consider the Hopf map $h : \S^3 \to \S^2$ given by $h(z,w) = (2\Real (z\bar{w}), \abs{z}^2 - \abs{w}^2)$. The preimage of any closed curve in $\S^2$ is called a Hopf torus~\cite{Pinkall1985}. In particular, the preimage of the circle with constant latitude $\delta-\pi/2$, i.e. $\abs{z}^2 - \abs{w}^2 = \sin(\delta - \pi/2)$, is a torus congruent to the image of $g$. For $\delta = \pi/2$, it is the preimage of the equator and thus the Clifford torus. As $\delta$ tends to $0$ or $\pi$, the circles of constant lattitude tend to the north and south poles and their corresponding Hopf tori tend to the two great circles in $\S^3$ that respectively lie over those poles.

% \begin{figure}
% \includegraphics[width=\textwidth]{graphics/genus0_linked.png}
% \caption{The image of two harmonic maps of the form~\eqref{eqn:genus zero simple map} with $\delta=\pi/4$ (black) and $3\pi/4$ (blue). The two images are interlinked, congruent and share a common circle $\Set{\exp(\omega \Sigma_2)}{\omega \in \R}$, drawn in red. In fact all tori of the form~\eqref{eqn:genus zero simple map} are tangent along this circle. In the limit as $\delta \to 0$ or $\delta \to \pi$ the image collapse to this circle. \\
% ~\\
% Every point of $\S^3$ except those on this circle belongs to exactly one image, so varying the parameter $\delta$ sweeps out all of $\S^3$.
% The image of a harmonic map with $\delta=\pi/2$ is the Clifford torus and it divides $\S^3$ into two congruent solid tori.\label{fig:genus0 linked}}
% \end{figure}

Thus far, the maps $g$ are maps from $\C$ to $SU(2)$. To descend to a torus $M$ it must be periodic with respect to a lattice $\Z\langle \tau_1,\tau_2\rangle$. We recall some formulae about matrix exponentials in $\SU(2)$. The standard inner product on the Lie algebra $\su_2$ is $\langle A, B \rangle_{\su_2} := -\tfrac{1}{2}\tr AB$. It follows that the norm $\norm{\cdot}$ of an $\su_2$ matrix is the square root of its determinant. This allows us to succinctly write 
\[
    \exp Z = I \cos \norm{Z} + \hat{Z} \sin \norm{Z},
    \;\;\;\;
    \text{where }\hat{Z} = \frac{1}{\norm{Z}} Z.
\labelthis{eqn:exp formula}
\]
An element $w\in\C$ is a period of $g$ exactly when
\begin{align*}
I = g(w) = I \cos&(4 w_R)\cos(4 w_I x)
- \hat{X}\sin(4 w_R)\cos(4 w_I x) \\
&+ \hat{Y}\cos(4 w_R)\sin(4 w_I x)
- \hat{X}\hat{Y}\sin(4w_R)\sin(4 w_I x),
\end{align*}
using~\eqref{eqn:exp formula}. The set $\{I,X,Y,XY\}$ is linearly independent, which forces 
% Squaring the coefficients of $\hat{Y}$ and $\hat{X}\hat{Y}$ and adding them together shows that $\sin^2(4w_Ix) = 0$. Doing likewise for the coefficients of $\hat{X}$ and $\hat{X}\hat{Y}$ shows that $\sin^2(4w_R) = 0$. Let 
$4w_R = \pi k$ and $4 w_I x = \pi l$ for integers $k$ and $l$. However, it also forces $(-1)^{k+l} = 1$.
Thus the lattice of the periods of $g$ is generated by
\[
\kappa_1 := \frac{\pi}{4}\bra{1 - \iu\frac{1}{x}},
\;\;\;\text{and}\;\;\;
\kappa_2 := -\frac{\pi}{4}\bra{1 + \iu\frac{1}{x}}.
\labelthis{eqn:def kappa12}
\]
The periods $\tau_1,\tau_2$ must be chosen from this lattice.
That is, there must be integers $n^1,m^1,n^2,m^2$
\[
\begin{pmatrix}
\tau_1 \\ \tau_2
\end{pmatrix}
= 
\begin{pmatrix}
n^1 & m^1 \\
n^2 & m^2
\end{pmatrix}
\begin{pmatrix}
\kappa_1 \\ \kappa_2
\end{pmatrix}
% \tau_1 = n^1 \kappa_1 + m^1 \kappa_2,\;\;
% \tau_2 = n^2 \kappa_1 + m^2 \kappa_2.
\labelthis{eqn:def tau12}
\]
These four integers may be interpreted as winding numbers of the map. The parallelogram spanned by $\kappa_1$ and $\kappa_2$ covers the image exactly once. Thus in~\eqref{eqn:def tau12} the integers $n^1$ and $m^1$ may be interpreted as how many times the loop $[0,\tau_1] \subset \C/\langle \tau_1,\tau_2 \rangle$ is wrapped around the image, and likewise for $n^2$ and $m^2$.
To have an image that is not a circle, the columns of the above matrix of integers should be linear independent. 

Varying $\delta$ changes the image of the map, however varying $x$ changes the domain torus. If we define $\tau = \tau_2/\tau_1$ to be the conformal type of the domain torus $M$, we arrive at the formula
\[
\tau
= \frac{(n^2 + m^2) + \iu x(n^2 - m^2)}
{(n^1+m^1) +\iu x(n^1-m^1)}.
\labelthis{eqn:conformal type}
\]
This shows that the conformal type of domain of the map $g$ depends only $x$ and four integers. One natural choice to make is $\tau_l=\kappa_l$ for $l=1,2$. That gives a conformal parameter of 
\[
\tau = \frac{1-x^2 + 2\iu x}{1+x^2}.
\]
As $x$ varies from $0$ to $\infty$, $\tau$ sweeps out the upper half of the unit circle. The range of the conformal parameter for any other map in our collection is the image of this semicircle under some element of $\SL_2\Q$, namely a circle centred on the real axis with rational endpoints (or a line in $\Q+\iu\R$).
A reason that the range of possible conformal parameters $\tau$ is so pleasing is that $\kappa_1$ and $\kappa_2$ always span a rhombus, which restricts the possible $M$ that admit a harmonic map of the form~\eqref{eqn:exp formula}.

We turn now to computing the spectral data associated to one of these harmonic maps $g: M \to \SU(2)$. Following~\cite{Hitchin1990}, we first give an $\S^1$ family of flat unitary connections. It is possible to analytically continue this representation to a $C^\times$ family of flat $\SL_2\C$ connections. From this family, we consider the holonomy representation of the fundamental group of $M$. As $M$ is a torus, the matrices of this representation commute and thus share eigenspaces. Generically, these eigenspaces will be two lines. It can be shown that the eigenlines coincide for only a finite number of points in $\C^\times$. If we consider the family of eigenlines as a subspace of $\C^\times \times \CP^1$, then there are only finitely many ramification points under projection to the first component. Further, it is possible to complete this to an algebraic curve in $\CP^1 \times \CP^1$ that double covers the first component. This is called the eigenline curve. For eigenline curves of genus 0, 1, and 2 this is exactly the spectral curve. For higher genera, one must take into account the possibility of eigenlines coinciding to higher order.

Concretely, for any harmonic map $g$ from $M$ to $\SU(2)$ the $\S^1$ family $d_\zeta$ of flat unitary connections is given by the formula
\[
d_\zeta := d_L + \frac{1}{2}g^{-1}dg + \zeta^{-1}\Phi - \zeta\Phi^*,
% = d + \bar{\kappa}^2 F\,d\bar{z} - \kappa^2 F^*\,dz + \zeta^{-1}F\,dz - \zeta F^*\,d\bar{z},
\labelthis{eqn:flat connection translation}
\]
where $\zeta\in\S^1$ is the parameter of the family, $d_L$ is the left invariant connection on $\SU(2)$, and $2\Phi$ is the $(1,0)$ part of $g^{-1}dg$. Because the form $g^{-1}dg$ is valued in $\su(2)$ we can write $g^{-1}dg = 2(\Phi - \Phi^*)$, and thus $d_\zeta = d_L + (1+\zeta^{-1})\Phi - (1+\zeta)\Phi^*$.
It is immediately clear how this extends analytically to a $\C^\times$ family. Taking $g$ to be a particular map of the form~\eqref{eqn:genus zero simple map}, we can compute that
\[
\Phi = \exp(-4w_I Y)(-Xdw -iYdw)\exp(4w_I Y).
\]
Next, we apply the gauge transformation $h = \exp(-4 w_I Y)$ to change the connection matrix $C$ to $h^{-1}dh + h^{-1}Ch$, which is constant in $w$.
% In our case, we may 
% \[
% d_\zeta = d -\frac{\iu}{2}(dw - d\bar{w}) + (1+\zeta^{-1})(-Xdw -iYdw) - (1+\zeta)(-X^*d\bar{w} -iY^*d\bar{w}).
% \]
As it is constant, for each $\zeta$ one can solve the parallel transport equation $d_\zeta V = 0$ explicitly by exponentiation.
The holonomy matrix for the loop from $0$ to $\tau_l$ is $H^l(\zeta) = \exp B^l(\zeta)$, for $B^l(\zeta)$ equal to 
\[
\zeta^{-1}(\tau_l + \bar{\tau_l}\zeta) 
\begin{pmatrix}
0 & - (1 + \iu xe^{\iu \delta}) + (-1 + \iu x e^{\iu \delta})\zeta \\
(1+\iu xe^{-\iu \delta}) + (1-\iu xe^{-\iu \delta})\zeta & 0
\end{pmatrix}.
\]
We note that our selection of the form of the vectors $X$ and $Y$ has ensured that the matrices are off-diagonal. Thus their eigenvalues and eigenspaces are simple to write down.

To find the spectral curve, we find the values of $\zeta$ for which the two eigenlines of $H^l(\zeta)$ coincide. The eigenspaces of $H^l(\zeta)$ and $B^l(\zeta)$ are the same, so we may do our computation with the latter. We have that $u(\zeta) = (u_1(\zeta),\; u_2(\zeta))^T$ is an eigenvector if and only if
\begin{multline*}
-\zeta^{-1}(\tau_l + \bar{\tau_l}\zeta)(1 - \iu x e^{\iu \delta})(\zeta-\alpha) u_2(\zeta)^2 \\
= \zeta^{-1}(\tau_l + \bar{\tau_l}\zeta)(1+\iu xe^{-\iu \delta}) (1-\bar{\alpha}\zeta) u_1(\zeta)^2.
\end{multline*}
where $\alpha$ is a point that is always inside the unit circle, given by
\[
\alpha = \frac{1+\iu x e^{\iu \delta}}{-1+\iu x e^{\iu \delta}}
= \frac{x e^{\iu \delta} - \iu}{x e^{\iu \delta} +\iu}.
\labelthis{eqn:def branch point genus zero}
\]
The points of the eigenline curve in $\CP^1 \times\CP^1$ are hence of the form
% \begin{multline*}
% \Bigg( \zeta, \Bigg[\pm \sqrt{-\zeta^{-1}(\tau_l + \bar{\tau_l}\zeta)(1 - \iu x e^{\iu \delta})(\zeta-\alpha)} \\
% : \sqrt{\zeta^{-1}(\tau_l + \bar{\tau_l}\zeta)(1+\iu xe^{-\iu \delta}) (1-\bar{\alpha}\zeta)} \Bigg] \Bigg)
% \end{multline*}
% \vspace{-0.2cm}
\[
% = 
\Bigg( \zeta, \Bigg[\pm \sqrt{-(1 - \iu x e^{\iu \delta})(\zeta-\alpha)} : \sqrt{(1+\iu xe^{-\iu \delta}) (1-\bar{\alpha}\zeta)} \Bigg] \Bigg).
\]
% \begin{align*}
% \Bigg( \zeta, \Big[\pm &\sqrt{-\zeta^{-1}(\tau_l + \bar{\tau_l}\zeta)(1 - \iu x e^{\iu \delta})(\zeta-\alpha)} \\
% &\hspace{3cm}: \sqrt{\zeta^{-1}(\tau_l + \bar{\tau_l}\zeta)(1+\iu xe^{-\iu \delta}) (1-\bar{\alpha}\zeta)} \Big] \Bigg) \\
% &= \Bigg( \zeta, \Big[\pm \sqrt{-(1 - \iu x e^{\iu \delta})(\zeta-\alpha)} : \sqrt{(1+\iu xe^{-\iu \delta}) (1-\bar{\alpha}\zeta)} \Big] \Bigg).
% \end{align*}
From this we can see where and to what order the eigenlines coincide. The plus-minus sign produces two distinct lines unless one of the components of the homogeneous coordinates has a root, with the order of coincidence the same as the order of the root. As $\iu x e^{\iu \delta}$ is always in the left half of the complex plane, $1 - \iu x e^{\iu \delta}$ and its conjugate $1 + \iu x e^{-\iu \delta}$ never vanish. Hence the eigenlines coincide only over $\alpha$ and $\cji{\alpha}$, and only to first order. The spectral curve is therefore $\eta^2 = (\zeta-\alpha)(1-\bar{\alpha}\zeta)$, a genus zero hyperelliptic curve without singularities.

Let us explore how variation of the parameter $\alpha$ may alter the properties of the harmonic map $g$, and provide some intuition about the limit as $\alpha$ approaches the unit circle. From~\eqref{eqn:def branch point genus zero}, we can see how the two continuous parameters $\delta$ and $x$ have been incorporated into the definition of $\alpha$. If we treat $xe^{\iu \delta}$ as a point in the upper half plane then~\eqref{eqn:def branch point genus zero} is the Cayley transform, a M\"obius transformation of the upper half plane to the unit disc. This shows that every $\alpha$ in the unit disc is obtained. 
One can write the inverse transformation as
\[
x e^{\iu \delta} = \iu \frac{1+\alpha}{1-\alpha}.
\]
Taking the magnitude of both sides shows that $x$ is constant along the arcs of circles centred on the real axis with radii such that the circle is perpendicular to the unit circle. 
If $x$ is constant, so too is $\tau$, and along this arc the corresponding family of harmonic maps have the same domain but varying images.
% as in Figures~\ref{fig:genus0 linked}.
In the limit as $\alpha$ approaches the unit circle for fixed $x$, the parameter $\delta$ tends to $0$ or $\pi$. Recall the $\delta$ is the angle between $X$ and $Y$, so that in the limit these vectors are parallel and the image of the harmonic map collapses into a great circle in $\SU(2)$.

Conversely if $\delta$ is fixed, say at $\delta=\pi/2$, but $x$ is allowed to vary then $\alpha = (x-1)/(x+1)$,
and takes values along the real axis. At the two extremes, when $\alpha$ tends to $-1$ or $1$, $x$ tends to $0$ or $\infty$ respectively. 
% Thus one period of the domain torus is dwarfing the other. 
% Another way to put this is that the limit of the lattice of periods will be only rank one, not rank two. 
However, throughout this deformation the image of the corresponding harmonic maps is fixed; it is only the domain that is changing. We could see this as one period of the domain becoming very large, but the derivative of the map in that direction becoming correspondingly small. Taken together, the result will be a map from the cylinder to a circle. This limiting process $\alpha \to \pm 1$ is not as well behaved as $\alpha \to \S^1\setminus\{\pm 1\}$ and we are not confident that this phenomena is replicated in higher spectral genera.
% TODO A cylinder only has one periodicity condition, one Sym point. Does \alpha tending to a Sym point blow it away, or make it the only important one?

Having found the spectral curve, we next compute the pair of differentials. The pair of differentials $\Theta^1$ and $\Theta^2$ arise as the derivatives of the logarithms of the eigenvalues $\mu^l$ of the holonomy matrices $H^l$ for $l=1,2$. Though the matrix $H^l$ has two eigenvalues, $\mu^l$ is a well defined function on the eigenline curve. Because $H^l$ is an $\SL_2\C$ matrix, its two eigenvalues are reciprocal. In this particular case, $(\mu^l)^{\pm 1} = \exp (\pm \nu^l)$, where $\nu^l$ is an eigenvalue of $B^l$, and $\Theta^l = d\log \mu^l = d\nu^l$. To compute $\nu^l$ we note that as $B$ is a traceless matrix,
\begin{align*}
(\nu^l)^2
= -\det B^l
&= -\zeta^{-2}(\tau_l + \bar{\tau_l}\zeta)^2 \abs{1- \iu xe^{\iu \delta}}^2 (\zeta-\alpha)(1-\bar{\alpha}\zeta).
\labelthis{eqn:eigenvalue}
\end{align*}
Therefore the differential $\Theta^l$ corresponding to the eigenvalues of $H^l(\zeta)$ is
\[
\Theta^l = d\,\log \mu^l = d\, \Big[ \zeta^{-1}(\tau_l + \bar{\tau_l}\zeta) \iu \abs{1 - \iu xe^{\iu \delta}} \eta \Big].
\]
We may use equation~\eqref{eqn:def branch point genus zero} to rewrite the terms with $x$ and $\delta$ to depend only on $\alpha$.
Recalling the definitions of $\kappa_1,\kappa_2$ from~\eqref{eqn:def kappa12}, we define scalars
\begin{align*}
r_l &:= \iu \kappa_l \abs{1-\iu x e^{\iu \delta}} = \frac{\pi}{2}\bra{ \frac{1}{\abs{1+\alpha}} + \iu \frac{1}{\abs{1-\alpha}} },
\,
\Psi^l := d\, \Big\{ \zeta^{-1}(r_l + \bar{r_l}\zeta) \eta \Big\}. 
\labelthis{eqn:genus zero differential basis}
\end{align*}
% Using these we may succinctly write $\Theta^l = n^l \Psi^1 + m^l \Psi^2$.
Hence the differentials lie in a lattice spanned by $\Psi^1, \Psi^2$ and we can identify it with the lattice of periods of the map $g$. This gives the interpretation that the pair of differentials in the spectral data determine the winding of the torus onto its image. This same interpretation holds for the general construction of a harmonic map from spectral data, where the domain of the map is constructed as the parallelogram spanned by the pair of differentials.
% \Psi is choosing a basis of the torus, \Theta is choosing the basis of the domain. It counts the winding with respect to these. It's all happening in the domain, but "a minimal region of the domain that covers the image". The map factors through the map spans by the \Psis. Have a look in lattice theory?
% Its the unique smallest/biggest lattice that supports a periodic map. The lattice IS unique, the basis is not

% TODO use the short exact to link differentials and lattice more directly

The final piece of the spectral data, though one that is less important in this paper, is the eigenline bundle $E$ on $\Sigma$. As $\Sigma$ is a sphere, up to isomorphism there is only one line bundle for each degree. By condition~\ref{P:quaternionic}, the line bundle $E$ must be degree $-1$, so there is a unique choice. Thus we have computed the spectral data for all maps $g$ of the form~\eqref{eqn:genus zero simple map} and shown that these are exactly the spectral data with a genus zero spectral curve.

We can also describe the space $\mathcal{M}_0$ of all harmonic maps of tori in $\mathbb{S}^3$ whose spectral curve is genus zero. As just noted, there is exactly one choice of line bundle, so the space is the same as the space of triples $(\Sigma,\Theta^1,\Theta^2)$. We have seen that the spectral curve $\Sigma$ is completely determined by its sole branch point $\alpha$ in the unit disc $D$. 
% And we have just seen that the differentials may be chosen rom 
%  rank two lattice. However by condition \ref{P:linear independence} they must be real linearly independent. Given a basis of the lattice of differentials, such as $\Psi^1$ and $\Psi^2$ in \eqref{eqn:genus zero differential basis} above, we may represent our choice of lattice points in terms of two pairs of integers. 
We can express the choice of differentials as an integer matrix equation
\begin{align*}
\begin{pmatrix}
\Theta^1 \\ \Theta^2
\end{pmatrix}
=
\begin{pmatrix}
n^1 & m^1 \\
n^2 & m^2
\end{pmatrix}
\begin{pmatrix}
\Psi^1 \\ \Psi^2
\end{pmatrix}.
\end{align*}
Linear independence, condition~\ref{P:linear independence}, is equivalent to this matrix having non-zero determinant. If we define $\Mat_2^*\Z = \Set{ M \in \Mat_2\Z }{\det M \neq 0}$. The moduli space $\mathcal{M}_0$ can be described succinctly as the product $D \times \Mat_2^*\Z$. 
One should bear in mind that this is not a canonical identification, as it is dependent on the choice of basis.

To finish this section, there is a formula to compute the energy of harmonic map from its spectral data, given in~\cite[Theorem 12.17]{Hitchin1990}, in terms of the spectral curve and differentials.
In particular, when the genus of the spectral curve is zero the differentials are entirely determined by the choice of the four integers $n^1,m^1,n^2, m^2$ and a point $\alpha$ in the unit disc, as in~\eqref{eqn:genus zero differential basis}, yielding
\[
E = \pi^2(1+\alpha\bar{\alpha})\frac{m^1 n^2 - n^1 m^2}{\abs{1-\alpha^2}}.
\]
The factor $1+\alpha\bar{\alpha}$ may be seen as a measure of how far the map is from being conformal; $\alpha=0$ is conformal and the energy increases as one moves away from this point. 
% Get rid of interpret, say it comes from the determinant and the size of the differentials
% One can interpret the fraction as giving the area of the domain of the map, this is the `area' spanned by the differentials. 
The fraction gives the `area' spanned by the differentials, which corresponds to the domain of the map.
% Tighten this up. \alpha=0 is conformal. Flag as a hand-wavy statement
As $\alpha$ tends to $\pm 1$, the energy grows without bound showing yet again that this limit is different to $\alpha$ tending to any other point on the unit circle.

\begin{figure}
\resizebox{0.5\textwidth}{!}{
\begin{tikzpicture}
\begin{axis}[
    title={$E(\alpha)$},
    xlabel=$\Real\,\alpha$, ylabel=$\Imag\,\alpha$,
]
\addplot3[
% 	surf,
% 	domain = -0.9:0.9,
% 	domain y=-1:1,
% ]
% 	{9*(1 + x^2 + y^2)/(sqrt((1-x^2+y^2)^2 + 4*x^2 *y^2))};
	surf,
    domain=0:0.9, %% radius
    samples=20,
    domain y=0:360, %% angle
    samples y=90,
    z buffer=sort,
    variable=\r,
    variable y=\t
]
(
{r*cos(t)},
{r*sin(t)},
{9*(1 + r^2)/(sqrt(1+r^2 - 2*r*cos(2*t)))}
);
\end{axis}
\end{tikzpicture}
}
\caption{A plot of the energy as a function over $\alpha \in D$.}
\end{figure}

%%%%%%%%%%%%%%%%%%%%%%%%%%%%%%%%%%%%%%%%%%%%%%%%%%%%%%%%%%%%%%%%%%%%%%%%
%%%%%%%%%%%%%%%%%%%%%%%%%%%%%%%%%%%%%%%%%%%%%%%%%%%%%%%%%%%%%%%%%%%%%%%%
%%%%%%%% GENUS ONE SECTION
%%%%%%%%%%%%%%%%%%%%%%%%%%%%%%%%%%%%%%%%%%%%%%%%%%%%%%%%%%%%%%%%%%%%%%%%
%%%%%%%%%%%%%%%%%%%%%%%%%%%%%%%%%%%%%%%%%%%%%%%%%%%%%%%%%%%%%%%%%%%%%%%%
\section{Spectral Curves of Genus One}\label{sec:Genus One}

The aim of this section is to describe the space $\mathcal{M}_1$ of harmonic maps of tori whose spectral curve has genus one. The main step is to describe the moduli space $\mathcal{S}_1$ of spectral curves of genus one, which we consider as a subspace of $\mathcal{A}/\mathbb{Z}_2$ where $\mathcal{A} := D\times D\setminus \Delta$ and $\Delta$ is the diagonal of the product. 
% We will construct the universal cover and recover $\mathcal{A}/\mathbb{Z}_2$ as the quotient of a certain subspace $\mathcal{\tilde{S}}$ of the by the group of covering transformations.
We will prove that the path connected components of $\mathcal{S}_1$ are indexed by two rational numbers $p > 0$ and $q$. For $p\neq 1$ the components are `helicoids' $(0,1)\times \R$, whereas for $p=1$ the components are annuli. We can then describe the bundle structure of $\mathcal{M}_1$ over $\mathcal{S}_1$.

There are essentially five stages in the proof of this result. 
For each real genus one curve there is a plane of real differentials meeting conditions~\ref{P:poles}--\ref{P:periods}. In Part~\ref{sub:Differentials} we find explicit formulae for these differentials with periods $0$ and $2\pi\iu$. 
The challenge is to determine whether it possible on a given curve to find a pair of differentials that additionally satisfies the closing conditions~\ref{P:closing}. This is taken up in Part~\ref{sub:Closing Conditions}.
We can choose an exact differential that satisfies the closing conditions exactly when the function $S$, defined by~\eqref{eqn:def_S}, is rationally valued.
Similarly, there is a non-exact differential on a curve that satisfies the closing condition exactly when the function $T$ is rational, where $T$ is defined by~\eqref{eqn:def_T}. This `function' $T$ is only well-defined locally, but the condition that $T$ is rational is well-defined so long as $S$ is rational. In Lemma~\ref{lem:closing_conds} we prove that a curve admits spectral data exactly when these two functions are rationally valued.

The third stage addresses this deficiency in $T$. In Part~\ref{sub:Universal Cover} we lift $T$ to a well-defined function $\tilde{T}$ on $\mathcal{\tilde{A}}$, the universal cover of $\mathcal{A}$. Let $\mathcal{\tilde{A}}(p)$ be the subspace where $S=p$ is fixed (we apply this notation generally). In Lemma~\ref{lem:T_graph} we show that for each $p$ there is a coordinate chart where $q=\tilde{T}$ is a coordinate. One has to consider $p \geq 1$ and $p \leq 1$ separately, but in both cases the coordinate charts have the form 
\[
\mathcal{\tilde{A}}(p) \diffeo \left\{ \bra{q,k,\tilde{X}} \in \R \times (0,1) \times \R \right\}.
\]
This lemma therefore proves that the level sets $S = p$ and $\tilde{T} = q$ in $\mathcal{\tilde{A}}(p)$, which we denote $\mathcal{\tilde{A}}(p,q)$, are contractible. Indeed, by our choice of coordinates they are coordinate planes. The decomposition of the preimage of the space of spectral curves into path connected components is
\[
\tilde{\pi}^{-1}(\mathcal{S}_1) = \coprod_{p \in \Q^+, q\in\Q} \mathcal{\tilde{A}}(p,q).
\]
The fourth stage, in Part~\ref{sub:Topology}, pushes this decomposition in the universal cover back down into the space of curves $\mathcal{A}/\mathbb{Z}_2$.
The group of covering transformations is proved to be $\langle \tilde{\lambda} \rangle$ and the action of $\tilde{\lambda}$ on $\mathcal{\tilde{A}}(p)$ is shown in Lemma~\ref{lem:T shift} to be
\[
\tilde{\lambda} : \bra{q,k,\tilde{X}} \mapsto \bra{q + (p-1),k,\tilde{X} + \pi}.
\]
Thus $\tilde{\lambda}$ maps $\mathcal{\tilde{A}}(p,q)$ to $\mathcal{\tilde{A}}(p,q + p-1)$. 
% If $p$ and $q$ are rational, so too is $q + (p-1)$ and hence the group $\mathcal{G} = \Z\langle \tilde{\lambda}\rangle$ of covering transformations restricts to give a group action on the preimage of the space of spectral curves.

Having described the path components of $\mathcal{S}_1$, in the final stage we can examine the fibre bundle structure of $\mathcal{M}_1$. Clearly it is trivial over contractible components of the base. However, over the components of $\mathcal{S}_1$ which are annuli, where $p=1$, the bundle is not trivial. We saw in the previous section that the fibre can be described in terms of matrices whose columns are independent elements of the $\Z^2$ lattice of differentials. The nontrivial structure of the bundle can be expressed in terms of a monodromy action on this lattice. The total space of the bundle $\mathcal{M}_1$ is shown to be a disjoint union of two dimensional strips.

Recall that $\mathcal{A}(p)/\mathbb{Z}_2$ denotes the subspace where $p$ is fixed. These coordinates now foliate of $\mathcal{A}/\mathbb{Z}_2$. Explicitly, for $p\neq 1$ in Theorem~\ref{thm:topology_curves} we have that
\[
\mathcal{A}/Z_2(p)
\diffeo \mathcal{\tilde{A}}(p)/\mathbb{Z}\langle\tilde{\lambda}\rangle
= \left\{ \bra{[q],k,\tilde{X}} \in \R/(p-1)\Z \times (0,1) \times \R \right\}.
\]
The subset where $[q] \in \Q/(p-1)\Z$ is $\mathcal{S}_1\cap\mathcal{A}(p)/\mathbb{Z}_2$.
This is a union of `helicoids'.
For $p=1$, by contrast, Theorem~\ref{thm:topology_curves_p1} shows that
\[
\mathcal{A}(1)/\mathbb{Z}_2
\diffeo \mathcal{\tilde{A}}(1)/\mathbb{Z}\langle\tilde{\lambda}\rangle
= \left\{ \bra{q,k,\left[\tilde{X}\right]} \in \R \times (0,1) \times \R/\pi\Z \right\}.
\]
This time $\mathcal{S}_1\cap\mathcal{A}(1)/\mathbb{Z}_2$ is the subset where $q\in\Q$ and is a union of annuli.

We summarise the spaces that we have introduced and their relationships to one another in the diagram below. The horizontal arrows represent covering maps, labelled with the group of covering transformations, whereas the vertical arrows represent inclusions. One starts out with the space $\mathcal{A}/\mathbb{Z}_2$ of real genus one curves in the top right corner. Inside here is the space of spectral curves $\mathcal{S}_1$. Along the top row we have that $\mathcal{A}/\mathbb{Z}_2$ is covered by the parameter space $\mathcal{A}$, which is in turn covered by the universal cover $\mathcal{\tilde{A}}$. Below that we have the subspaces $\mathcal{\tilde{A}}(p)$ on which the function $S$ has the value $p\in\Q^+$. On the bottom line we have the statement that the preimage of the space of spectral curves is a union of level sets $\mathcal{\tilde{A}}(p,q)$, the subsets of $\mathcal{\tilde{A}}$ on which $S = p$ and $\tilde{T} = q$. 

\[
\begin{diagram}
    \mathcal{\tilde{A}} &\rOnto^{2\pi \Z}&  \mathcal{A}:= \mathcal D\times D\setminus \Delta  &\rOnto^{\Z_2}&  \mathcal{A}/\mathbb{Z}_2 \\
    \uInto &&  &&  \\
    %%%%%%%%%%%%%%%%%%%%%%%%%%%%%%%%%%%%%
    \coprod_{p\in\Q^+} \mathcal{\tilde{A}}(p)  &&    && \uInto \\
    \uInto  &&  &&  \\
    %%%%%%%%%%%%%%%%%%%%%%%%%%%%%%%%%%%%%
    \tilde{\pi}^{-1}(\mathcal{S}_1) = \coprod_{p\in\Q^+, q\in\Q} \mathcal{\tilde{A}}(p,q)  && \rOnto && \mathcal{S}
\end{diagram}
\]

\subsection{Explicit Formulae for the Differentials}\label{sub:Differentials}
For subsequent calculations, it will be necessary to have explicit formulae for the differentials on a curve that have integral periods in Jacobi normal form. We will define two differentials $\Theta^E$ and $\Theta^P$ such that any differential $\Theta$ that meets conditions~\ref{P:poles}--\ref{P:periods} may be written as a linear combination of these two.

With the standard scaling, our curves of genus one may be written as
\[
\eta^2 = P(\zeta) := (\zeta-\alpha)(1-\bar{\alpha}\zeta)(\zeta-\beta)(1-\bar{\beta}\zeta).
\labelthis{eqn:genus one curve}
\]
One can equivalently identify this curve with its branch points $\alpha,\beta$ inside the unit circle. However $(\alpha, \beta)$ and $(\beta,\alpha)$ in $\mathcal{A} := (D \times D) \setminus \Delta$ (where $\Delta$ is the diagonal) both yield the same curve. Thus each curve is a point in $\mathcal{A}/\mathbb{Z}_2$. 

On the other hand, every elliptic curve may be transformed into the Jacobi normal form $w^2 = (1-z^2)(1-k^2z^2)$, where $k$ is a complex number called the elliptic modulus. 
Choosing the M\"obius transformation $f$ is essentially choosing a correspondence for the roots, and the value of $k$ may be determined from the cross ratio $[\alpha,\cji{\alpha};\beta,\cji{\beta}]$.
There are twenty-four possible choices, but we choose the correspondence
\[
  \begin{array}{ l||c|c|c|c}
    \zeta & \alpha & \cji{\alpha} & \beta & \cji{\beta} \\
    \hline
    z = f(\zeta) & 1 & -1 & k^{-1} & -k^{-1}
  \end{array}.
  \labelthis{eqn:choice f}
\]
This correspondence has three properties that distinguish it from the others. By convention, $k \in (0,1)$, which rules out sixteen of the choices. Second, consider the behaviour of the curve as $k\to 1$. The Jacobi form of the curve develops two nodes at $z=\pm 1$. This corresponds to forming nodes $\alpha=\beta$ and $\cji{\alpha} = \cji{\beta}$. This choice also takes the interior of the unit disc to the right half plane, which is our preferred convention. There is only one other correspondence that has these three properties, with the difference being which root, $\alpha$ or $\beta$, is mapped to $1$. In order to make this choice consistently as the branch points vary one must work with $\mathcal{A}$ and not $\mathcal{A}/\mathbb{Z}_2$.

% Considering the cross ratio of the branch points
% \[
% [\alpha,\cji{\alpha};\beta,\cji{\beta}] = \frac{\abs{\alpha-\beta}^2}{\abs{1-\bar{\alpha}\beta}^2}. \labelthis{eqn:roots_cross_ratio}
% \]
% This is a real quantity, and so the four roots lie on a circle (or a line). Thus the roots of the Jacobi form do also, which forces $k$ to be real. Any transformation between the curves must take branch points to branch points, 
% Choosing

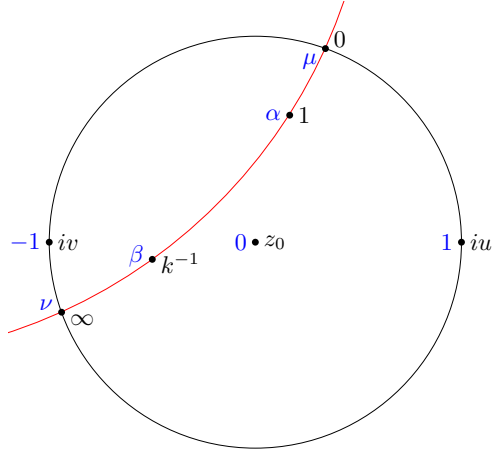
\begin{figure}
    \resizebox{0.56\textwidth}{!}{
        \begin{tikzpicture}
            \clip (-3.6,-3.5) rectangle (3.6,3.5);
            
            \draw (0,0) circle [radius=3];
            \draw[color=red] (-6,6) circle [radius=7.7];
            
            \fill (0,0) circle (0.05) 
                node[left,color=blue]{$0$}
                node[right,color=black]{$z_0$};
            \fill (3,0) circle (0.05) 
                node[left,color=blue]{$1$}
                node[right,color=black]{$\iu u$};
            \fill (-3,0) circle (0.05) 
                node[left,color=blue]{$-1$}
                node[right,color=black]{$\iu v$};
            
            % \fill (2.7,-3) circle (0.05) node[left, color=red]{$\infty$} node[right,color=black]{$-\bar{z}_0$};
            
            \fill (1.02,2.82) circle (0.05)
                node[below=5pt, left, color=blue]{$\mu$}
                node[above=4pt, right,color=black]{$0$};
            \fill (0.5,1.85) circle (0.05)
                node[left, color=blue]{$\alpha$}
                node[right,color=black]{$1$};
            \fill (-1.5,-0.25) circle (0.05)
                node[above=1.5pt, left, color=blue]{$\beta$}
                node[below=1.5pt, right,color=black]{$k^{-1}$};
            \fill (-2.82,-1.02) circle (0.05)
                node[above=3pt,left, color=blue]{$\nu$}
                node[below=3pt,right,color=black]{$\infty$};
            
            \end{tikzpicture}
    }      
        \caption{The $\zeta$-plane, with points marked in blue. The black labels are their images under $f$. The red line is the circle through the branch points.\label{fig:zeta plane}}
    \end{figure}

From our correspondence, equating the cross ratios $[\alpha,\cji{\alpha};\beta,\cji{\beta}]$ and\\ $[1,-1;k^{-1},-k^{-1}]$ gives
\[
k = \frac{\abs{1-\bar{\alpha}\beta}-\abs{\alpha-\beta}}{\abs{1-\bar{\alpha}\beta}+\abs{\alpha-\beta}},
\labelthis{eqn:def_k}
\]
and the map $f$ can be computed from the relation
\[
[\alpha,\cji{\alpha};\beta,\zeta] = [1,-1;k^{-1},f(\zeta)].
\]
Instead of solving this relation for $f$ immediately, we first look at the geometry of the transformation.
The image under $f$ of unit circle in the $\zeta$-plane is the imaginary axis in the $z$-plane. The involution $\rho(\zeta)$ fixes the unit circle, and exchanges the pairs of branch points $\alpha,\cji{\alpha}$ and $\beta,\cji{\beta}$. The corresponding antiholomorphic involution $\tilde{\rho}(z) := -\bar{z}$ in the $z$-plane exchanges $1,-1$ and $k^{-1},-k^{-1}$. The four branch points of the spectral curve lie on a circle, which is mapped to the real axis by $f$. Let the two points at the intersection of this circle with the unit circle be $\mu$ and $\nu$, with $\mu$ lying between $\alpha$ and $\cji{\alpha}$ and $\nu$ lying between $\beta$ and $\cji{\beta}$ (see Figure~\ref{fig:zeta plane}). Hence $f(\mu)$ and $f(\nu)$ must lie on the intersection of the real and imaginary axes, with $f(\mu) = 0$ and $f(\nu) = \infty$.

A M\"obius transformation, such as $f$, is determined up to scaling by the points it sends to $0$ and $\infty$, in this case $\mu$ and $\nu$. One other point is therefore needed to determine this scaling. We write $z_0 := f(0)$. Using the reality structure $\tilde{\rho}(z)$, we have $f(\infty) = -\bar{z_0}$. These points allow us to write concise formulae for $f$ and $f^{-1}$
\begin{align}
z = f(\zeta) &= -\bar{z}_0 \frac{\zeta - \mu}{\zeta - \nu},
\label{eqn:f} \\
\zeta = f^{-1}(z) &= \nu \frac{z - z_0}{z + \bar{z_0}}.
\label{eqn:f_inv}
\end{align}
Because of the holomorphic involution $\sigma: \eta\to-\eta$, equations~\eqref{eqn:f} and~\eqref{eqn:f_inv} almost but not quite specify a relation between $\eta$ and $w$: there is a free sign choice to make. 
% On the curve $\Sigma = \{ (\zeta,\eta) | \eta^2 = P(\zeta) \}$ there are two disjoint circles in $\Sigma$ lying over the unit circle in $\CP^1$.
% At a point $(\zeta,\eta)$ over the unit circle in $\Sigma$, we have that the value of $\eta$ is $\pm \zeta\abs{\zeta-\alpha}\abs{\zeta-\beta}$.
% Thus there is a notion of the `positive' unit circle, the one on which $\eta$ is positive over $\zeta=1$.
We define the function $\eta^+(\zeta)$ on the unit circle to be $\zeta\abs{\zeta-\alpha}\abs{\zeta-\beta}$ and $w^+(\iu u) = + \sqrt{1+u^2}\sqrt{1+k^2u^2}$.
% Under the transformation $f$ the unit circle is mapped to the imaginary axis. At the points over the imaginary axis in the Jacobi elliptic curve, those points $(z,w)$ for which $z=\iu u$, we have that $w = \pm \sqrt{1+u^2}\sqrt{1+k^2u^2}$. Again, it is possible to make a consistent choice of sign along these two disjoint circles. 
For the sake of being concrete, we choose the transformation between elliptic curves that maps $\eta^+$ to $w^+$.
% the positive unit circle to points over the imaginary axis where $w$ is positive. 
In a slight abuse of notation, we shall also use $f$ to denote the map between elliptic curves.

Having found the map $f$ that transforms a genus one curve into Jacobi normal form, we may return our attention to the differentials that satisfy conditions~\ref{P:poles}--\ref{P:periods}.
There is one differential that obviously meets these requirements, namely
\[
\Theta^E := \iu\; d\left( \frac{\eta}{\zeta} \right).
% = \iu\left[ -\alpha\beta + \frac{1}{2}\left(\alpha(1+\abs{\beta}^2) + \beta(1+\abs{\alpha}^2)\right)\zeta - \frac{1}{2}\left(\bar{\alpha}(1+\abs{\beta}^2) + \beta(1+\abs{\alpha}^2)\right)\zeta^3 + \bar{\alpha}\bar{\beta}\zeta^4 \right]\frac{d\zeta}{\zeta^2\eta}.
\labelthis{eqn:thetaE}
\]
The superscript $E$ is a mnemonic for exact. We therefore seek non-exact differentials that also meet these requirements.
% It is standard to refer to differentials with double poles and no residues as differentials of the second kind. Condition~\ref{P:poles} may therefore be rephrased as the differential must be of the second kind and have poles at the points of $\Sigma$ over $\zeta=0$ and $\infty$. 
The standard Jacobi differential of the first kind is simply the holomorphic differential $\omega := \tfrac{dz}{w}$.
The standard Jacobi differential of the second kind is defined to be
$e := (1-k^2 z^2) \tfrac{dz}{w}$.
Every differential of the second kind may be written as the linear combination of $e$, the holomorphic differential $\omega$, and an exact differential~\cite[Art. 167]{Hancock1910}.
The differential $e$ is real with respect to $\tilde{\rho}(z) = -\bar{z}$, but it has a pole at $z=\infty$ (i.e. $\zeta=\nu$), which not allowed under condition~\ref{P:poles}. This pole can be moved to $\zeta=\infty$ by adding an exact differential. Hence we define 
\[
\varepsilon := e + d\left[ \frac{w}{z + \bar{z}_0} \right] + \Real z_0\; d\left[ \frac{w}{(z-z_0)(z + \bar{z}_0)} \right]
= e + d\left[ \frac{(z - \iu\,\Imag z_0)w}{(z-z_0)(z + \bar{z}_0)} \right].
\labelthis{eqn:def_varepsilon}
\]
There is essentially a `freedom' to add any real multiple of $\Theta^E$ to this differential and it will still satisfy~\ref{P:poles}--\ref{P:reality}. This particular definition of $\varepsilon$ however has a principal part perpendicular to the principal part of $\Theta^E$, which introduces a symmetry that we will employ later.

We now compute the periods of the differentials $\omega$ and $\varepsilon$. 
The standard choice of branch cuts in the $z$-plane is $[1,k^{-1}]$ and $[-k^{-1},-1]$. The standard basis of homology is a clockwise loop around $-1$ and $1$, the $A$ period, and an anticlockwise loop around $1$ and $k^{-1}$, the $B$ period. This $B$ period is homologous to a the imaginary axis traversed top to bottom. 
The $A$ periods of the standard elliptic differentials define functions of the elliptic modulus called the complete elliptic integrals of the first and second kind, denoted $K$ and $E$ respectively.
\begin{align*}
    \int_{A} \omega &=: 4K(k), &
    \int_{A} e &=: 4E(k).
\end{align*}
Generally we omit the elliptic modulus $k$.
The expression $\sqrt{1-k^2}$ is known as the complementary modulus and by convention we write $K' = K(\sqrt{1-k^2})$ and $E' = E(\sqrt{1-k^2})$, where the prime represents complement not derivative. 
Its relevance is due to Jacobi's imaginary transformation~\cite[p. 34]{Jacobi1829}, which shows the $B$ periods are
\begin{align*}
\int_{B} \omega &= 2\iu K', &
\int_{B} e &= 2\iu(K'-E').
\end{align*}

Note that $\varepsilon$ is the sum of $e$ and an exact differential and so has the same periods as $e$. We define $\Theta^P = 2E\omega - 2K\varepsilon$, or if we unwind the definition~\eqref{eqn:def_varepsilon} of $\varepsilon$,
\[
\Theta^P = 2E \omega - 2Ke - 2K d\left[ \frac{(z-\iu\,\Imag z_0)w}{(z-z_0)(z + \bar{z}_0)} \right].
\labelthis{eqn:thetaP}
\]
The superscript $P$ is a mnemonic for period, because this differential has an $A$ period of $0$ and a $B$ period of $2\pi\iu$. This follows from the above standard periods and from Legendre's relation $K'E + KE' - KK' = \frac{\pi}{2}$,
% TODO reference    
which is essentially a reciprocity law~\cite{Griffiths1994}.
% TODO reference
The above equation for $\Theta^P$ shows a nice division, with the first two terms providing the desired periods and the last term giving the required poles. Though this definition may seem arbitrary, we can characterise the differential $\Theta^P$ in the following natural way.

\begin{lem}\label{lem:theta2_characterisation}
The differential $\Theta^P$ is the unique real differential on the curve $\Sigma$ with double poles and no residues over $\zeta=0$ and $\infty$, with periods $0$ and $2\pi\iu$ on $A$ and $B$ respectively, and with its principal part over $\zeta=0$ satisfying
\[
\pp \Theta^P \in \iu \R \pp \Theta^E.
\]
\begin{proof}
% All differentials on a genus one curve $\Sigma$ that have (at worst) a double pole over $\zeta=0$ and $\zeta=\infty$ may be written the form
% \[
% \Theta = b(\zeta)\frac{d\zeta}{\zeta^2\eta},
% \]
% for a polynomial $b(\zeta)$ of degree $4$. If we write $b(\zeta) = b_0 + \cdots + b_4 \zeta^4$, Condition~\ref{P:reality} forces $b_i = b_{4-i}$. Lastly, if we write out the equation of $\Sigma$ as $\eta^2 = P(\zeta) = P_0 \zeta + \cdots + P_4 \zeta^4$, $\Theta$ has no residues exactly when $P_1b_0 - 2P_0b_1 = 0$. 

% We first verify that $\Theta^P$ has such properties, then verify uniqueness. The only property not yet demonstrated is the third one, concerning the principal part. We note that $\zeta=0$ corresponds to $z=z_0$, so for some real scalar $r$
% \[
% \pp \Theta^E
% = \iu \pp d\bra{\frac{\eta}{\zeta}}
% = \iu \pp d \bra{ \frac{rw}{(z-z_0)(z + \bar{z}_0)} }
% = -\iu r \frac{w(z_0)}{z_0 + \bar{z}_0}\frac{dz}{(z-z_0)^2}.
% \]
% And on other side we have
% \begin{align*}
% \pp \Theta^P
% &= - 2K \pp d\left[ \frac{(z-\iu\,\Imag z_0)w}{(z-z_0)(z + \bar{z}_0)} \right] \\
% &= + 2K \frac{(z_0-\iu\,\Imag z_0)w(z_0)}{z_0 + \bar{z}_0} \frac{dz}{(z-z_0)^2} \\
% &= 2K (\Real z_0) \frac{w(z_0)}{z_0 + \bar{z}_0} \frac{dz}{(z-z_0)^2}.
% \end{align*}
Suppose that $\Theta$ is a differential with the properties granted in the lemma. Then $\Theta-\Theta^P$ is exact, real, and has double poles with no residues over $\zeta = 0,\infty$. Choose $s\in\R$ such that $\Theta-\Theta^P - s\Theta^E$ has no pole over $\zeta = 0$. By reality it is holomorphic. The only exact holomorphic differential is zero, and hence $\Theta = \Theta^P + s \Theta^E$.
As taking principal part is a linear operation, the final condition forces $s=0$.
\end{proof}
\end{lem}

Using the same proof as the above lemma, any differential satisfying~\ref{P:poles}--\ref{P:periods} may be written in the form $a \Theta^E + n \Theta^P$ for some $a\in\R$ and $n\in\Z$,

\subsection{The Closing Conditions}\label{sub:Closing Conditions}
We have seen that every curve $\Sigma\in\mathcal{A}/\mathbb{Z}_2$ admits differentials that meet Conditions~\ref{P:poles}--\ref{P:periods}, but it is not possible to find differentials on every curve that further satisfy the closing conditions~\ref{P:closing}. We will define functions on $\mathcal{A}$ such that the conditions are satisfiable for a curve exactly when the value of the functions are both rational.

The closing conditions on the spectral data ensure that it corresponds to a harmonic map of the torus, rather than a merely a harmonic section of some bundle on the torus.
Let us recall~\ref{P:closing}, supposing that $\Theta$ satisfies~\ref{P:poles}--\ref{P:periods}. The closing condition at $\zeta=1$ if
\[
\int_{\gamma^{+}} \Theta \in 2\pi\iu \Z,
\]
where $\gamma^+$ is a path that begins at $(1,-\eta^+(1))$ and ends at $(1,\eta^+(1))$. Although the value of the integral is dependent on the path, due to~\ref{P:periods} the condition is well-defined. Likewise, the closing condition at $\zeta=-1$ is defined by integrating over a path $\gamma^-$ from $(-1,-\eta^+(-1))$ to $(-1,\eta^+(-1))$. 

First we can reduce to a special case. Suppose that we have a triple of spectral data $(\Sigma,\Theta^1,\Theta^2)$ such that the differentials $\Theta^1$ and $\Theta^2$ have imaginary periods $2\pi\iu l_1$ and $2\pi\iu l_2$ respectively. Let $l>0$ be the greatest common denominator of $l_1$ and $l_2$, and by B\'ezout's identity let $x$ and $y$ be integers that satisfy $xl_1 + yl_2 = l$. Then consider the differentials $\Psi^E$ and $\Psi^P$ defined by
\[
\vt{\Psi^E}{\Psi^P} =
\begin{pmatrix}
\tfrac{l_2}{l}    &   -\tfrac{l_1}{l} \\
x                       &   y
\end{pmatrix}
\vt{\Theta^1}{\Theta^2}.
\]
The imaginary periods of this new pair of differentials are $0$ and $2\pi\iu l$ respectively, and they also meet the closing conditions. Further, the two differentials are linearly dependent exactly when $\Psi^E$ is zero, compare to~\ref{P:linear independence}. Hence, a curve admits spectral data if and only if it admits spectral data where the differentials have imaginary periods $0$ and $2\pi\iu l$,.

For an exact differential, necessarily $a\Theta^E$ for some $a\in\R$, the particular path of integration is irrelevant and, using~\eqref{eqn:thetaE}, the closing conditions are
\[
\int_{\gamma^{+}} a\Theta^E 
= a\iu \left. d\bra{\frac{\eta}{\zeta}} \right|_{(1, -\eta^+(1))}^{(1, \eta^+(1))} 
= 2\iu \eta^+(1)a 
% = 2i \abs{1-\alpha}\abs{1-\beta}a 
\in 2\pi\iu \Z,
\]
and $-2\iu \eta^+(-1)a \in 2\pi\iu\Z$.
% \[
% \int_{\gamma^{-}} a\Theta^E = -2\iu \eta^+(-1)a = 2\iu \abs{1+\alpha}\abs{1+\beta}a \in 2\pi\iu \Z.
% \]
There is a common solution for $a$ if and only if
\[
S(\alpha,\beta) := \frac{2\iu \eta^+(1)}{-2\iu \eta^+(-1)} = \frac{\abs{1-\alpha}\abs{1-\beta}}{\abs{1+\alpha}\abs{1+\beta}} \in \Q^+.
\labelthis{eqn:def_S}
\]
% The condition above, $S\in\Q^+$, is a necessary condition for a spectral curve to admit spectral data. 

To find a second function, concerning the existence of a differential with imaginary period $2\pi\iu l$, we follow a similar line of reasoning. For some real number $b$, we may write $\Psi^P = b \Theta^E + l \Theta^P$. Fix two paths $\gamma^+, \gamma^-$. The two closing conditions applied to $\Psi^P$ are then
\begin{equation}
\left.
\begin{aligned}
2\iu \eta^+(1) b + l\int_{\gamma^+} \Theta^P &=: 2\pi\iu \Gamma^+ \in 2\pi\iu \Z, \\
-2\iu \eta^+(-1) b + l\int_{\gamma^-} \Theta^P &=: 2\pi\iu \Gamma^- \in 2\pi\iu \Z.
\end{aligned}
\quad
\right\}
\label{eqn:period closing cond}
\end{equation}
Elimination of $b$ yields the condition for a common solution to exist. This condition can be written as
\[
2\pi\iu T(\alpha,\beta) := S(\alpha,\beta) \int_{\gamma^-} \Theta^P - \int_{\gamma^+} \Theta^P
% = 2\pi \iu \frac{S(\alpha,\beta) \Gamma^- - \Gamma^+}{l}
\in 2\pi\iu \Q.
\labelthis{eqn:def_T}
\]
% using the definition of $S(\alpha,\beta)$ to substitute for $-\eta^+(1)/\eta^+(-1)$. 
% As was shown at the beginning of this section, whether the integral of a differential over $\gamma^+$ or $\gamma^-$ lies in $2\pi\iu \Z$ is independent on the particular path chosen between the marked points. In the same manner, 
This time, the value of the function is dependent on the particular paths chosen; 
$T(\alpha,\beta)$ is a multi-valued function on $\mathcal{A}$ defined up the addition of an element from $\Z\langle 1,S(\alpha,\beta)\rangle$.
However, if $S(\alpha,\beta)$ is rational then the condition $T(\alpha,\beta)\in\Q$ is independent of the choice of paths.

\begin{lem}\label{lem:closing_conds}
A curve admits spectral data if and only if $S\in\Q^+$ and $T\in\Q$. 
Hence, the space of genus one spectral curves is the projection of the following set to $\mathcal{A}/\mathbb{Z}_2$:
\[
\{(\alpha,\beta) \in \mathcal{A} \mid S(\alpha,\beta) \in \Q^+, T(\alpha,\beta) \in \Q\}
\]

\begin{proof}
From the above discussion, the rationality of $S$ and $T$ are a necessary conditions.
For the converse, we will explicitly construct a pair of differentials. From the previous section, we can always find a pair of non-zero differentials with properties~\ref{P:poles}--\ref{P:periods}, namely $\Theta^E$ and $\Theta^P$, but these do not necessarily satisfy~\ref{P:closing}. We will find $a,b \in \R$ and $l\in \N$ such that $\Psi^E = a\Theta^E$ and $\Psi^P = a\Theta^E + l\Theta^P$ that further satisfy~\ref{P:closing} and~\ref{P:linear independence}, which are all the conditions necessary to be spectral data. In fact, we will find the minimal such differential, in the sense that all other differentials satisfying~\ref{P:poles}--\ref{P:closing} are an integer combination of $\Psi^E$ and $\Psi^P$.

Fix $\gamma^+$ and $\gamma^-$ so that $T$ is well-defined and let $S = n/m$ and $T = n'/m'$ for coprime integers.
First we will produce the exact differential $\Psi^E$. 
% The rationality of $S$ directly ensures the consistent solution of an $a$.
We define
\[
a := \frac{2\pi\iu n}{2\iu \eta^+(1)} = \frac{2\pi\iu n}{-2\iu \eta^+(-1)}\frac{1}{S} = \frac{2\pi\iu m}{-2\iu \eta^+(-1)},
\]
and observe that for $\Psi^E := a\Theta^E$,
\begin{align*}
\int_{\gamma^{+}} \Psi^E &= 2\iu \eta^+(1)a = 2\pi \iu n, 
&
\int_{\gamma^{-}} \Psi^E &= -2i \eta^+(-1)a = 2\pi \iu m.
\end{align*}
This demonstrates the existence of an exact differential $\Psi^E$ that satisfies the closing conditions~\ref{P:closing}. 
% By assumption $n$ is not zero, so neither is $\Psi^E$.
Any other exact differential $\Theta$ that meets the closing conditions must also be a multiple of $\Theta^E$. If its integrals of $\gamma^+,\gamma^-$ are $\Gamma^+,\Gamma^-$ respectively, then 
\[
\frac{n}{m} 
= \frac{\int_{\gamma^+} \Theta^E}{\int_{\gamma^-} \Theta^E} 
= \frac{\int_{\gamma^+} \Theta}{\int_{\gamma^-} \Theta} 
= \frac{\Gamma^+}{\Gamma^-}.
\]
Since $n/m$ is a simplified fraction we must have that $\Gamma^+ = cn$ for an integer $c$. Hence $\Theta = ca\Theta^E = c \Psi^E$.

To find a non-exact differential $\Psi^P$ that also satisfies the closing conditions, we must solve~\eqref{eqn:period closing cond} for $b$ and $l$ such that $\Gamma^+$ and $\Gamma^-$ are integral. Eliminating $b$ from~\eqref{eqn:period closing cond} gives
\begin{align*}
% \frac{1}{2\iu \eta^+(1)} \bra{2\pi\iu \Gamma^+ - l\int_{\gamma^+} \Theta^P}
% &= \frac{1}{-2\iu \eta^+(-1)} \bra{2\pi\iu \Gamma^- - l\int_{\gamma^-} \Theta^P} \\
% 2\pi\iu \Gamma^+ - l\int_{\gamma^+} \Theta^P
% &= \frac{2\iu \eta^+(1)}{-2\iu \eta^+(-1)} \bra{2\pi\iu \Gamma^- - l\int_{\gamma^-} \Theta^P} \\
2\pi\iu \bra{ S \Gamma^- - \Gamma^+ }
&= l \bra{S \int_{\gamma^-} \Theta^P - \int_{\gamma^+} \Theta^P} \\
% \frac{ n \Gamma^- - m \Gamma^+ }{m}
% &= l \frac{n'}{m'} \\
m'(n \Gamma^- - m \Gamma^+)
&= lmn'.
\labelthis{eqn:period diff existence}
\end{align*}
By B\'ezout's identity, since $m$ and $n$ are coprime there exists an integral solution for $\Gamma^+$ and $\Gamma^-$ exactly when $m'$ divides the right hand side. The smallest $l$ for which this holds is $m' / \gcd(m',mn')$. We may define $b$ from either equation of~\eqref{eqn:period closing cond} and set $\Psi^P = b\Theta^E + l\Theta^P$. Explicitly 
\[
b := \frac{1}{2\iu \eta^+(1)}\bra{ 2\pi\iu \frac{mn'}{\gcd(m',mn')}y - l \int_{\gamma^+} \Theta^P }.
\labelthis{eqn:def Psi coeff}
\]
As $\Psi^E$ is exact and $\Psi^P$ is non-exact, the two differentials are real linearly independent (Condition~\ref{P:linear independence}).

Finally if $\Theta$ is any differential that satisfies the closing conditions, we may write it as $\Theta = b'\Theta^E + l'\Theta^P$ for $b'\in\R$ and $l'\in \Z$. Its imaginary period is $2\pi\iu l'$ and let its integrals over $\gamma^+$ and $\gamma^-$ be $2\pi\iu \Gamma^+$ and $2\pi\iu \Gamma^-$. These coefficients must satisfy Equation~\eqref{eqn:period diff existence} and hence $l$ divides $l'$. Because the differential $\Theta - \frac{l'}{l}\Psi^P$ meets conditions~\ref{P:poles}--\ref{P:closing} and is exact, it must be an integer multiple of $\Psi^E$ by the first part of this proof. Therefore $\Theta$ is an integer combination of $\Psi^E$ and $\Psi^P$. 
\end{proof}
\end{lem}

% For higher genus spectral curves, not considered in this paper, the closing conditions are difficult to work with, harder than even the period conditions. In the genus one case the integrals of $\Theta^E$ lead to an algebraic expression, as we have just seen in~\eqref{eqn:def_S}, but the integrals of $\Theta^P$ will lead to a transcendental conditions involving incomplete elliptic integrals.

Immediately before this lemma we observed that $T(\alpha,\beta)$ is a multi-valued function on $\mathcal{A}$. 
On each curve of $\mathcal{A}\setminus\{\nu = \pm 1\}$ we make the choice of paths $\gamma_0^+$ and $\gamma_0^-$ (note the bold symbols) as per Figure~\ref{fig:gamma paths}. This effectively chooses a principal branch $T_0$ of $T$. Note that these paths depend on the order of the branch points $(\alpha,\beta) \in \mathcal{A}$, since $\nu$ is defined to be the point on the unit circle that lies between $\beta$ and $\cji{\beta}$. We now develop a formula for $T_0$ suitable for calculations.

% \labelpara{para:principal paths}
\begin{figure}
\resizebox{0.6\textwidth}{!}{
    \begin{tikzpicture}

\clip (-3.5,-3.5) rectangle (3.5,3.5);

\draw[->] (-3.3,0)--(3.3,0) node[above left]{$\Real z$};
\draw[->] (0,-3.3)--(0,3.3) node[below right]{$\Imag z$};

\fill (0,2) circle (0.05) node[left]{$f(1)$};
\fill (0,-2.4) circle (0.05) node[left]{$f(-1)$};

\fill (2.1,0) circle (0.05) node[below=2.1pt,color=black]{$1$};
\fill (-2.1,0) circle (0.05) node[below=2.1pt,color=black]{$-1$};

\draw[color=red,dashed,>->]
    (-0.1,2.1)
    --(-0.1,1.1)
    node[left, color=red]{$\gamma_0^+$}
    --(-0.1,-0.1)
    --(2.3,-0.1)
    arc (-90:90:0.1);
\draw[color=red,->]
    (2.3,0.1)
    --(0.1,0.1)
    --(0.1,2.1);

\draw[color=blue,dashed,>->]
    (-0.2,-2.1)
    --(-0.2,-1.1)
    node[left, color=blue]{$\gamma_0^-$}
    --(-0.2,0.2)
    --(2.3,0.2)
    arc (90:-90:0.2);
\draw[color=blue,->]
    (2.3,-0.2)
    --(0.1,-0.2)
    --(0.1,-2.1);
\end{tikzpicture}
}
\caption{For a curve $\Sigma$ corresponding to $(\alpha,\beta)\in\mathcal{A}\setminus\{\nu = \pm 1\}$, consider the coordinate change given by $f$. Recall $f(\alpha)=1$. The principal path $\gamma_0^+$, in red, starts at $f(1)$ on the lower sheet, traverses around $1$ (crossing the branch cut) and then returns to $f(1)$ on the upper sheet, without crossing $\iu\infty = f(\nu)$. Likewise, the principal path $\gamma_0^-$ is marked in blue.\label{fig:gamma paths}}
\end{figure}

The incomplete Legendre elliptic integrals $F(z;k)$ and $E(z;k)$ are defined respectively to be the integrals of $\omega$ and $e$ from $0$ to $z$ on the $w^+$ sheet. Recalling the definition of $\Theta^P$ from Equation~\eqref{eqn:thetaP}, we have
\begin{align*}
\int_{\gamma_0^+} \Theta^P
&= 4 E F(f(1);k) - 4 K E(f(1);k) - 4K \frac{(f(1)-\iu\Imag z_0)w(f(1))}{(f(1)-z_0)(f(1) + \bar{z}_0)}, \\
\int_{\gamma_0^-} \Theta^P
&= 4 E F(f(-1);k) - 4 K E(f(-1);k) 
% \\
% &\qquad\qquad\qquad\qquad 
- 4K \frac{(f(-1)-\iu\Imag z_0)w(f(-1))}{(f(-1)-z_0)(f(-1) + \bar{z}_0)}.\labelthis{eqn:thetaPformulae}   
\end{align*}
% These two formulae may be substituted into~\eqref{eqn:def_T_0} to compute $T_0$. In the next section however, we will make a change of coordinates that simplifies these formulae.

% Previously, the comment was made that the particular algorithm to chose a path is not valid when $\nu = \pm 1$. Indeed, the result of this can be seen directly in the formulae we have derived. When $\nu$ takes either of these values, then one of $f(1)$ or $f(-1)$ will be infinite. We also note that these integrals are purely imaginary, as we expected on theoretical grounds, because $f(1)$ and $f(-1)$ are purely imaginary, $F(z;k)$ and $E(z ;k)$ take the imaginary axis to itself and
% \[
% (f(1)-z_0)(f(1) + \bar{z}_0) = -(f(1)-z_0)(\overline{f(1) - z_0}) = - \abs{f(1) - z_0}^2.
% \]
% The function $T_0$ is therefore real valued.

% \section{Coordinates for \texorpdfstring{$\mathcal{A}$}{A}}

It is clear from these formulae that the points $f(1)$ and $f(-1)$ play an important role in determining the value of $T_0$. It therefore prudent to adopt a parametrisation of $\mathcal{A}$ that is better suited to calculations than $(\alpha,\beta)$. We take three coordinates to be $k$, $\iu u = f(1)$ and $\iu v = f(-1)$. For the final coordinate we shall take $p=S(\alpha,\beta)$ itself, as then we can enforce the condition $S\in\Q^+$ simply by holding this coordinate constant. 
It is also necessary to exclude the plane where $u=v$. By the definition just given, equality is impossible since $f$ is invertible. In a limiting sense, these points correspond to the diagonal $\{\alpha=\beta\} \subset D\times D$. 
A technicality is that these coordinates $(p,k,u,v)$ only cover part of $\mathcal{A}$ since, for example, there are points of $\mathcal{A}$ where $\iu u = f(1)$ is infinite. This can be addressed by considering the imaginary axis together with infinity as $\iu\R\cup\{\infty\}=\RP^1$ and using $u^{-1}$ and $v^{-1}$ as coordinates where necessary. With this understanding, we claim
\[
\mathcal{A} = \left\{ (p,k,u,v) \in \R^+ \times (0,1) \times \RP^1\times \RP^1 \mid u\neq v\right\}.
\]
The transition functions between the $(\alpha, \beta)$ and $(p,k,u,v)$ coordinates have in part already been described; $p = S(\alpha,\beta)$ is given by Equation~\eqref{eqn:def_S} and $k$ by Equation~\eqref{eqn:def_k}. The map $f$ is also determined by $(\alpha,\beta)$ and as above $u = -\iu f(1)$ and $v = -\iu f(-1)$. In the reverse direction, if we can describe the map $f$ in terms of  $(p,k,u,v)$ and then $\alpha = f^{-1}(1)$ and $\beta = f^{-1}(k^{-1})$. As a M\"obius transformation is determined up to a scalar by the points sent to $0$ and $\infty$ we shall find an expression for $z_0 = f(0)$ in terms of the new coordinates. Note the following identity using cross ratios:
\[
\abs{\frac{\alpha-1}{\alpha+1}}
% = \abs{\frac{\alpha-1}{\alpha+1}} \abs{\frac{0+1}{0-1}}
= \abs{ \cross{\alpha}{0}{1}{-1} }
= \abs{ \cross{1}{z_0}{\iu u}{\iu v} }
= \abs{\frac{1-\iu u}{1 - \iu v}} \abs{\frac{z_0 - \iu v}{z_0 - \iu u}}
\]
Applying the same identity gives a similar formula for $\beta$. Together they show that
\[
p = S(\alpha,\beta)
= \abs{\frac{\alpha-1}{\alpha+1}} \abs{\frac{\beta-1}{\beta+1}}
= \abs{\frac{1-\iu u}{1 - \iu v}}\abs{\frac{1 - k\iu u}{1 - k\iu v}} \abs{\frac{z_0 - \iu v}{z_0 - \iu u}}^2,
\]
which, holding $(p,k,u,v)$ constant, forces $z_0$ to lie on a circle centred on the imaginary axis. On the other hand, in the $\zeta$-plane the points $-1$, $0$, and $1$ all lie on a straight line that is invariant under the real involution $\rho$. The corresponding statement in the $z$-plane is that $\iu v$, $z_0$, and $\iu u$ all lie on a circle that is symmetric under reflection in the imaginary axis. Thus we have determined two circles that $z_0$ lies on. As these two circles are both centred on the imaginary axis, they intersect in two points: $z_0$ and $-\bar{z_0}$. We may solve for $z_0$, giving
\[
z_0 = \frac{\sqrt{pw(\iu u)w(\iu v)}}{pw(\iu v) + w(\iu u)} \abs{u-v} + \iu\frac{puw(\iu v) + vw(\iu u)}{pw(\iu v) + w(\iu u)},
\labelthis{eqn:formula xy}
\]
where $w(z)^2 = (1-z^2)(1-k^2z^2)$. Finally the scaling of $f^{-1}(z)$ can be determined from the fact that $f^{-1}(\iu u) = 1$, hence 
\[
f^{-1}(z)
=  \frac{\iu u + \bar{z_0}}{\iu u - z_0} \frac{z-z_0}{z + \bar{z_0}}
% =  -\frac{\iu v + \bar{z_0}}{\iu v - z_0} \frac{z-z_0}{z + \bar{z_0}}.
\]
Since the formulae for $z_0$ and $f^{-1}$, and hence for $\alpha$ and $\beta$, are smooth in $(p,k,u,v)$ we have established this is a valid coordinate change on $\mathcal{A}$. From equations~\eqref{eqn:def_T} and~\eqref{eqn:thetaPformulae}, we now can write $T_0$ in terms of $(p,k,u,v)$ as
\begin{align*}
2\pi\iu T_0(p,k,u,v)
% &= p\left[ 4E F(\iu v;k) - 4K E(\iu v;k) - 4\iu K \frac{w(\iu v)}{u-v} \right]
% - \left[ 4E F(\iu u;k) - 4K E(\iu u;k) + 4\iu K \frac{w(\iu u)}{u-v} \right] \\
&= 4p \left[ E F(\iu v;k) - K E(\iu v;k) \right] - 4\left[ E F(\iu u;k) - K E(\iu u;k) \right] \\
&\qquad\qquad - 4\iu K \frac{p w(\iu v) + w(\iu u)}{u-v}.
\labelthis{eqn:Teqn}
\end{align*}

\subsection{The Universal Cover}\label{sub:Universal Cover}

The function $T_0$ that we have just constructed is apt for computation, but because $T$ is a multi-valued function any results on its level sets in $\mathcal{A}$ will be limited to local results. The way forward is to transition to the universal cover $\mathcal{\tilde{A}}$ of $\mathcal{A}$, which is covered by coordinates $(p,k,\tilde{u},\tilde{v})$. The universal cover allows us to pull back the multi-valued function $T$ to a single valued function $\tilde{T}$. This global function $\tilde{T}$ will allow us to gain global results about the space of spectral curves.

\begin{lem}\label{lem:mathcal tilde C}
The universal cover of $\mathcal{A}$ is
\[
\mathcal{\tilde{A}} =
\{(p, k,\tilde{u},\tilde{v}) \in \R^+\times(0,1)\times\R\times\R \mid  \tilde{u} < \tilde{v} < \tilde{u} + 2\pi \},
\]
with the projection map $\tilde{\pi} : \mathcal{\tilde{A}} \to \mathcal{A}$ given by
\begin{align*}
p = p,\quad 
k = k,\quad
u = \tan \frac{\tilde{u}}{2},\;
u^{-1} = \cot \frac{\tilde{u}}{2},\quad
v = \tan \frac{\tilde{v}}{2},\;
v^{-1} = \cot \frac{\tilde{v}}{2}.
\end{align*}

\begin{proof}
The coordinates $(p,k,u,v)$ embed $\mathcal{A}$ into $\R^+\times(0,1)\times \RP^1 \times \RP^1$. Its universal cover is $\R^+\times(0,1)\times\R^2$, with the covering map given above using the $2\pi$-periodic half-$\tan$ mapping of $\R$ to $\RP^1$. Next recall that $\mathcal{A}$ is the complement of the hyperplane $u-v = 0$.
Considering the $\RP^1\times\RP^1$ as a torus, the line $u=v$ is where the toroidal and poloidal angles are equal; removing this line leaves an annulus.
When the hyperplane $u-v=0$ is pulled back to the universal cover it becomes a collection of hyperplanes $\tilde{u}-\tilde{v} \in 2\pi\Z$. Thus we may take a simply connected region of the complement $\tilde{u}-\tilde{v} \not\in 2\pi\Z$ to cover $\mathcal{A}$ and this is $\mathcal{\tilde{A}}$ above.
\end{proof}
\end{lem}

We now lift and extend $T_0$ to a single valued function $\tilde{T}$ on $\mathcal{\tilde{A}}$. Recall that $T_0$ was defined on $\mathcal{A}\setminus\{\nu = \pm 1\}$, which is open dense in $\mathcal{A}$. The main difficulty is extending the incomplete integrals $F(z;k)$ and $E(z;k)$. The function $E(z;k)$ has a pole at infinity, so we work instead with $E(z;k) - kz$ which is well behaved. Both of these integrals, for $z\in \iu\R$, take purely imaginary values,
\[
F(\iu x; k)
= \iu \int_0^x \frac{dt}{\sqrt{(1+t^2)(1+k^2 t^2)}},
\]
and
\[
E(\iu x; k) - k\iu x
= \iu \int_0^x \frac{1-k^2}{\sqrt{1+t^2}\left( \sqrt{1+k^2t^2} + k\sqrt{1+t^2} \right)} \;dt.
\]
The imaginary parts are odd and increasing functions of $x\in\R$, so there is no way to extend them to single valued functions on $\RP^1$. Following the universal cover defined above, let $\tilde{x}\in\R$ and $x = \tan\frac{\tilde{x}}{2}$. We now explicitly construct $\tilde{F}(\tilde{x};k)$ and $\tilde{E}(\tilde{x};k)$ to be the analytic continuations of these integrals. Take $\tilde{x}$ in the range $(-\pi,\pi)$ and pull back through the covering map to define
\begin{align*}
\tilde{F}(\tilde{x};k)
&:= \int_0^{2 \atan \tilde{x}} \frac{dt}{\sqrt{(1+t^2)(1+k^2t^2)}}
% &= \int_0^{\tilde{x}} \frac{1}{\sqrt{(1+(\tan s/2)^2)(1+k^2(\tan s/2)^2)}} \times \frac{1}{2}\sec^2\frac{s}{2}\;ds \\
= \frac{1}{2} \int_0^{\tilde{x}} \frac{ds}{\sqrt{\cos^2 s/2 + k^2 \sin^2 s/2}} \;, \\
\tilde{E}(\tilde{x}; k)
% &= \int_0^{2\atan \tilde{x}} \frac{\sqrt{1+k^2t^2} - k\sqrt{1+t^2}}{\sqrt{1+t^2}}\;dt \\
% &= \frac{1}{2} \int_0^{\tilde{x}} \frac{\sqrt{\cos^2 s/2 + k^2\sin^2 s/2} - k}{ \cos^2 s/2 } \;ds \\
% &= \frac{1}{2} \int_0^{\tilde{x}} \frac{\bra{\cos^2 s/2 + k^2\sin^2 s/2} - k^2}{ \cos^2 s/2 } \frac{1}{\sqrt{\cos^2 s/2 + k^2\sin^2 s/2} + k}\;ds \\
&:= \frac{1}{2} (1-k^2) \int_0^{\tilde{x}} \frac{ds}{\sqrt{\cos^2 s/2 + k^2\sin^2 s/2} + k}\;.
\end{align*}
% Manipulating the integrand into a form that is plainly analytic demonstrates that $\tilde{E}$ is as well. 

\begin{figure}
\includegraphics[height=0.3\textheight]{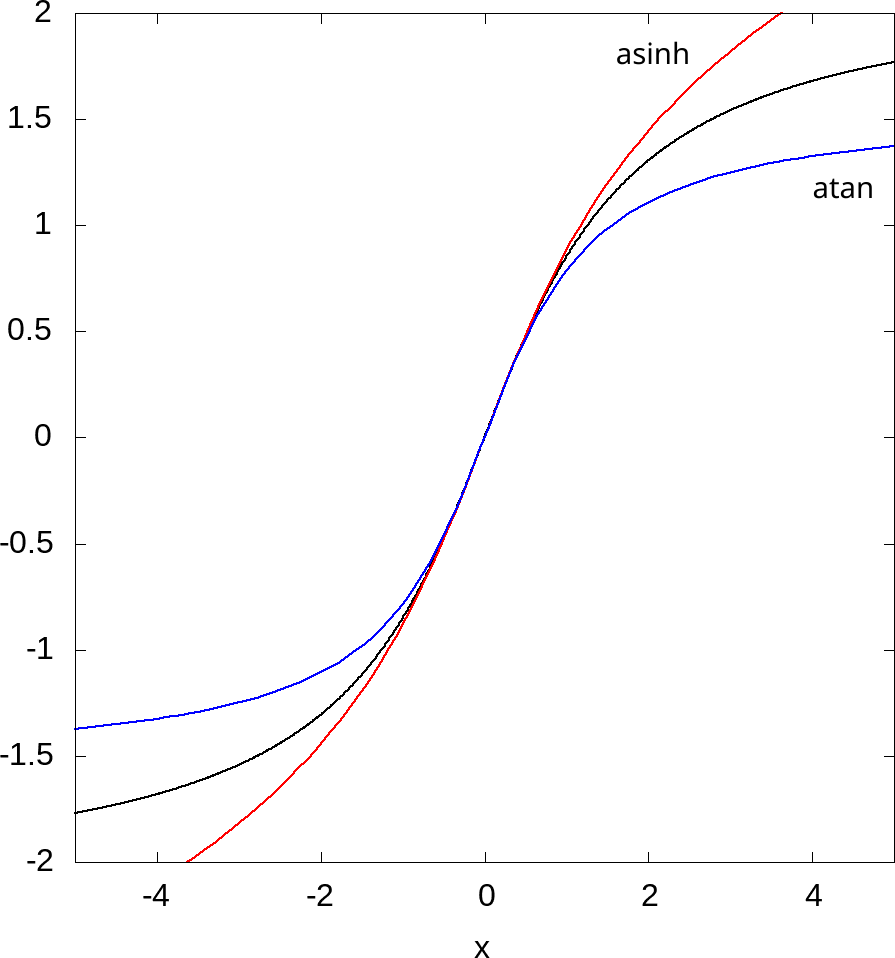}
\caption{Plot of $\Imag F(\iu x;0.5)$ (black), $\Imag F(\iu x;1) = \atan(x)$ (blue), and $\Imag F(\iu x;0)=\asinh(x)$ (red). For each value of $k$ the function $\Imag F(\iu x;k)$ is bounded, but as $k\to 0$ it tends to the unbounded $\asinh x$.}
\end{figure}
\begin{figure}
\includegraphics[height=0.3\textheight]{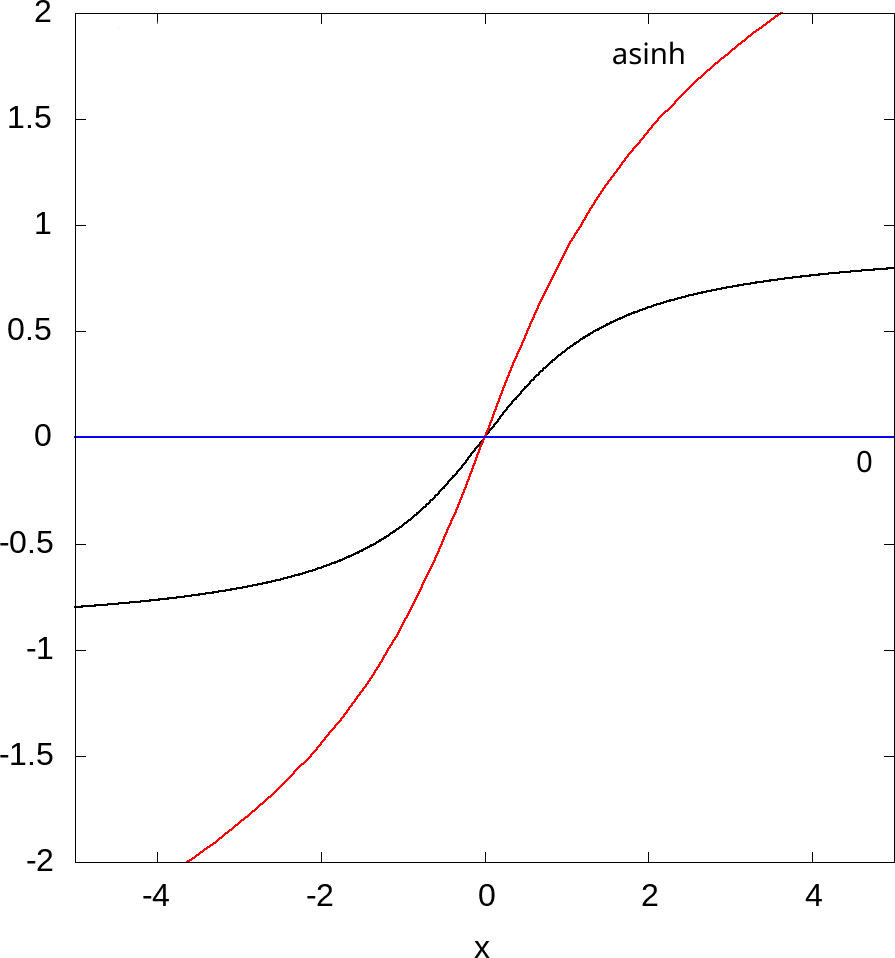}
\caption{Plot of $\Imag E(\iu x; 0.5) - 0.5x$ (black), $\Imag E(\iu x; 1) - x = 0$ (blue) and $\Imag E(\iu x; 0)=\asinh(x)$ (red). Like $F$, it is bounded for each value of $k$ but unbounded in the limit $k \to 0$.}
\end{figure}

If we rewrite $T_0$ in the following way,
\begin{align*}
2\pi T_0(p,k,u,v) 
= 4p& \left[ E \Imag F(\iu v; k) - K \Imag \left(E(\iu v; k) - k\iu v\right) \right] \\
 -4& \left[ E \Imag F(\iu u; k) - K \Imag \left(E(\iu u; k) - k\iu u\right) \right] \\
 -4& K \left[p\left(\frac{w(\iu v)}{u-v} + kv \right) + \left(\frac{w(\iu u)}{u-v} - ku \right) \right],
\labelthis{eqn:T0 rewrite}
\end{align*}
then the first two brackets can be lifted to the universal cover by substitution of $\tilde{F}$ and $\tilde{E}$. The third bracket of~\eqref{eqn:T0 rewrite} is analytic on $\mathcal{A}$ and so lifts to the universal cover directly via pulling back. 
To actually compute the value of the function $\tilde{T}$ note that the integrands of $\tilde{F}$ and $\tilde{E}$ are $2\pi$ periodic and the periods of the integrals are $2K'$ and $2(E'-K')$ respectively, from the $B$-periods of the standard elliptic integrals. By Legendre's relation
\[
E \tilde F(\tilde{x} + 2\pi; k) - K \tilde E(\tilde{x}+2\pi; k)
= E \tilde F(\tilde{x}; k) - K \tilde E(\tilde{x}; k) + \pi.
\]
Define the winding number $\Wind : \R \to \Z$ of a number $\tilde{x}$ to be the integer $\Wind(\tilde{x})$ such that $-\pi < \tilde{x} - 2\pi\Wind(\tilde{x}) < \pi$  and write $\tilde{x} = 2\pi \Wind(\tilde{x}) + x$ for $x = \tan \tilde{x}/2$. Then
% \[
% \tilde{F}(\tilde{x};k)
% % = \frac{1}{2} \bra{ \int_0^{2\pi \Wind(\tilde{x})} + \int_{2\pi \Wind(\tilde{x})}^{2\pi \Wind(\tilde{x}) +y}} \frac{1}{\sqrt{\cos^2 s/2 + k^2 \sin^2 s/2}} \;ds
% % &= \frac{1}{2} \bra{ n\int_0^{2\pi} + \int_{0}^{y}} \frac{1}{\sqrt{\cos^2 s/2 + k^2 \sin^2 s/2}} \;ds \\
% % &= \bra{ n\oint_{\RInf} + \int_{0}^{x}} \frac{dt}{\sqrt{(1+t^2)(1+k^2t^2)}} \\
% = 2\Wind(\tilde{x})K' + F_0(x;k),
% \qquad
% \tilde{E}(\tilde{x}; k) = 2n(K'-E') + E_0(x;k),
% \labelthis{eqn:tildeF_periodic}
% \]
% and
% \begin{align*}
% 2\pi \tilde{T}(p,k,\tilde{u},\tilde{v})
% %%%%%%%%%%%%%%%%%%%%%%%%
% % &= 4p \left[ E (2K'W(\tilde{v})+F_0(v)) - K (2(K'-E')W(\tilde{v})+E_0(v)) \right] \\
% % &\quad - 4 \left[ E (2K'W(\tilde{u})+F_0(u)) - K (2(K'-E')W(\tilde{u})+E_0(u)) \right]
% % - 4 K \left[p\left\{\frac{w(\iu v)}{u-v} + kv \right\}
% % + \left\{\frac{w(\iu u)}{u-v} - ku \right\} \right] \\
% %%%%%%%%%%%%%%%%%%%%%%%%
% % &= 2\pi T_0(p,k,u,v) + 4p \left[ 2EK' - 2K(K'-E') \right]\Wind(\tilde{v}) \\
% % &\qquad - 4 \left[ 2EK' - 2K(K'-E') \right]\Wind(\tilde{u}) \\
% %%%%%%%%%%%%%%%%%%%%%%%%
% \tilde{T}(p,k,\tilde{u},\tilde{v})
% &= T_0(p,k,u,v) + 2\left[ p\Wind(\tilde{v}) - \Wind(\tilde{u}) \right],
% \labelthis{eqn:tilde T computable}
% \end{align*}
\begin{align*}
2\pi \tilde{T}(p,k,\tilde{u},\tilde{v}) 
&:= 4p \left[ E \tilde{F}(\iu v; k) - K \tilde{E}(\iu v; k)\right]
-4 \left[ E \tilde{F}(\iu u; k) - K \tilde{E}(\iu u; k)\right] \\
&\qquad- 4 K \left[p\left(\frac{w(\iu v)}{u-v} + kv \right) + \left(\frac{w(\iu u)}{u-v} - ku \right) \right],\\
\tilde{T}(p,k,\tilde{u},\tilde{v})
&= T_0(p,k,u,v) + 2\left[ p\Wind(\tilde{v}) - \Wind(\tilde{u}) \right],
\labelthis{eqn:tilde T computable}
\end{align*}
Thus to do computations with the function $\tilde{T}$, for the most part one can continue to work with the function $T_0$ downstairs on $\mathcal{A}$ using~\eqref{eqn:Teqn} and keep track of the winding numbers. 

We use this strategy to prove three facts about $\tilde{T}$. First if $p$ is rational, $T\circ\tilde{\pi} \in \Q$ if and only if $\tilde{T} \in \Q$. This is immediate from the above formula. The second fact is Lemma~\ref{lem:range_T} concerning the range of $\tilde{T}$. The third is that it has non-zero derivatives under certain hypotheses, Lemma~\ref{lem:dTdu_nonzero}.

\begin{lem}\label{lem:range_T}
The range of \,$\tilde{T}$ on $\mathcal{\tilde{A}}$ is $\R$.

\begin{proof}
Note that the lemma is obvious then $p\neq 1$ by increasing $\tilde{u}$ and $\tilde{v}$ in increments of $2\pi$ in Equation~\eqref{eqn:tilde T computable}. However we shall prove below that the range of $T_0$ on $\mathcal{\tilde{A}} \setminus \{\nu = \pm 1\}$ is $\R$. It will follow then that the range of $\tilde{T}$ is also $\R$ for any $p$. Fix $p$ and $k$. From the imaginary periods, $\abs{F(\iu x;k)}$ is bounded by $K'$ and $\abs{E(x;k) - kx}$ is bounded by $K'-E'$. Thus there is some constant, dependent on both $p$ and $k$ but independent of $u$ and $v$ that bounds the first two terms of $T_0$ in~\eqref{eqn:T0 rewrite}.
% \begin{gather*}
% - 4 K \left[\frac{pw(\iu v) + w(\iu u)}{u-v} + k(pv-u) \right] - C \\
% \leq
% 2\pi T_0(p,k,u,v)
% \leq \\
% - 4 K \left[\frac{pw(\iu v) + w(\iu u)}{u-v} + k(pv-u) \right] + C,
% \end{gather*}
Therefore it is sufficient to show that for any fixed $p$ and $k$ the final terms have range equal to the real line. But this is easy to show. Consider the limit as $u \to v^\pm$,
\[
\lim_{u \to v^\pm} \frac{pw(\iu v) + w(\iu u)}{u-v} + k(pv-u)
= k(p-1)v + (p+1) w(\iu v)\lim_{u \to v^\pm} \frac{1}{u-v} = \pm\infty.
\]
% From the other side,
% \[
% \lim_{u \to v^-} \frac{pw(\iu v) + w(\iu u)}{u-v} + k(pv-u)
% = k(p-1)v + (p+1) w(\iu v)\lim_{u \to v^-} \frac{1}{u-v} = -\infty.
% \]
By continuity $T_0$ obtains every value.
\end{proof}
\end{lem}

We now compute the derivatives of $\tilde{T}$ and provide conditions under which they are non-zero. First note that for $\tilde{u} \not\in \pi + 2\pi\Z$,
\[
\Partial{\tilde{T}}{\tilde{u}} 
= \frac{1}{2}\sec^2\left(\frac{\tilde{u}}{2}\right) \Partial{T_0}{u}
= \frac{1}{2}(1+u^2) \Partial{T_0}{u}.
\] 
This is because (fixing $p$) the difference between the multivalued function $T$ and $T_0$ is locally constant, which is removed by differentiation.
From the explicit formula in~\eqref{eqn:Teqn} and recalling $w(\iu u) = \sqrt{(1+u^2)(1+k^2 u^2)}$, the $u$-derivative of $T_0$ is
\begin{align*}
\frac{\pi}{2}&w(\iu u)(u-v)^2\Partial{T_0}{u} \\
% &= -\frac{E}{w(\iu u)} + \frac{pK w(\iu v)}{(u-v)^2} + \frac{K}{w(\iu u)(u-v)^2}\left[1 + u^2 - uv + v^2 + k^2 uv + k^2 u^2v^2 \right]\\
&= -(u-v)^2 E + pKw(\iu u)w(\iu v) + K\left[ 1 + u^2 - uv + k^2 uv + v^2 + k^2 u^2 v^2 \right]
\labelthis{eqn:dTdu}
\end{align*}
We would like to say that the derivative never vanishes but that is not true. However it is true that it never vanishes when $p \geq 1$. We will address the $p < 1$ case in~\eqref{eqn:T_0 symmetry} in a different way.

\begin{lem}\label{lem:dTdu_nonzero}
If $p \geq 1$ then the $\tilde{u}$-derivative of $\tilde{T}$ is never zero.

\begin{proof}
We first show that the derivative of $T_0$ is non-zero where it is defined. Apply the crude estimate that $K>E>1$,
which follows form $K(0)=E(0)=\pi/2$ and that they are increasing and decreasing functions respectively, to obtain
\begin{align*}
\frac{\pi}{2}&w(\iu u)(u-v)^2\Partial{T_0}{u}\\
% &= -(u-v)^2 E + pKw(\iu u)w(\iu v) + K\left[ 1 + u^2 - uv + k^2 uv + v^2 + k^2 u^2 v^2 \right] \\
&> -(u-v)^2 K + Kw(\iu u)w(\iu v) + K\left[ 1 + u^2 - uv + k^2 uv + v^2 + k^2 u^2 v^2 \right] \\
&= K \left[ w(\iu u)w(\iu v) + 1 + (1 + k^2) uv + k^2 u^2 v^2 \right]
\end{align*}
This formula is almost sufficient to show the derivative is positive. The only term that could be negative is the one featuring $uv$. However, a lower bound for the square root terms is
\[
w(\iu u) = \sqrt{1 + (1+k^2)u^2 + k^2u^4} > \sqrt{(1+k^2)}\abs{u},
\]
so applying this to both $w(\iu u)$ and $w(\iu v)$ we have
\begin{align*}
\frac{\pi}{2}w(\iu u)(u-v)^2\Partial{T_0}{u}
% &> K \left[ (1+k^2)\abs{uv} + 1 + (1 + k^2) uv + k^2 u^2 v^2 \right] \\
&\geq K \left[ 1 + k^2 u^2 v^2 \right] >0.
\end{align*}
This shows that if $p\geq 1$ and $\tilde{u} \not\in \pi + 2\pi\Z$, the $\tilde{u}$ derivative of $\tilde{T}$ never vanishes. It remains to check it does not vanish when $\tilde{u} \in \pi + 2\pi\Z$. As $\tilde{T}$ is an analytic function its derivatives are continuous and so we may compute their value at these points by taking a limit.
\begin{align*}
\lim_{\tilde{u}\to \pi + 2\pi\Z} \Partial{\tilde{T}}{\tilde{u}}
% &=\frac{1}{\pi} \lim_{u \to \infty} \frac{1 + u^2}{w(\iu u)(u-v)^2} L\bra{ p,k,u,v } \\
% &=\frac{1}{\pi} \lim_{u \to \infty} \frac{(1 + u^2)u^2}{w(\iu u)(u-v)^2} \times u^{-2}L\bra{ p,k,u,v } \\
&=\frac{1}{\pi} \frac{1}{k} \left(-E + pkKw(\iu v) + K\left[ 1 + k^2v^2 \right]\right),
\end{align*}
The proof of the positivity of the above is immediate using $K>E>1$.
\end{proof}
\end{lem}

% This result is of interest because it shows that if $p \leq 1$, then the $v$ and $v'$ derivatives of $T$ are non-zero:
% \begin{align*}
% \frac{\pi}{2}\Partial{T}{v} &= \frac{-p}{w(\iu v)(u-v)^2} U\bra{ \frac{1}{p},k,u,v }, \\
% \frac{\pi}{2}\Partial{T}{v'} &= \frac{p}{w'(v')(uv'-1)^2} V\bra{ \frac{1}{p},k,u,v' }.
% \end{align*}
% And if $p \geq 1$, then the $u$ and $u'$ derivatives of $T$ are non-zero:
% \begin{align*}
% \frac{\pi}{2}\Partial{T}{u} &= \frac{1}{w(\iu u)(u-v)^2} U\bra{ p,k,u,v }, \\
% \frac{\pi}{2}\Partial{T}{u'} &= \frac{-1}{w'(u')(1-u'v)^2} V\bra{ p,k,v,u' }.
% \end{align*}

\subsection{The Topology of the Moduli Space}\label{sub:Topology}

In this final part of Section~\ref{sec:Genus One} we describe the space of genus one spectral curves $\mathcal{S}_1$. 
More precisely,
in Lemma~\ref{lem:T_graph} we demonstrate that for any values $p \in \R^+$, $q\in\R$ the level set defined by $p = S$ and $q = \tilde{T}$ is a graph over two coordinates. 
% This follows by an application of the implicit function theorem applied to $\tilde{T}$. Applying the implicit function theorem requires the derivatives of $\tilde{T}$ (given by Equation~\eqref{eqn:Teqn}) and showing that they are non-vanishing (Lemma~\ref{lem:deriv no zeroes}).
Hence each level set is diffeomorphic to $(0,1)\times\R$.
% Moreover, that lemma shows that for fixed $p$ these strips $\mathcal{\tilde{A}}(p,q)$ foliate $\mathcal{\tilde{A}}(p)$, the subspace of $\mathcal{\tilde{A}}$ on which $p$ is a fixed constant.
It follows from the observation that $\tilde{T}$ is rational exactly when $T$ is that the preimage of $\mathcal{S}_1$ in the universal cover $\mathcal{\tilde{A}}$ may be written as the disjoint union of the level sets.
% \[
% \mathcal{\tilde{S}} = \coprod_{p\in \Q^+,\; q\in \Q} (0,1)\times\R.
% \]
We recover $\mathcal{S}_1$ by taking the quotient of the level sets by the group of covering transformations of $\mathcal{\tilde{A}}$ over $\mathcal{A}/\mathbb{Z}_2$. This culminates in Theorems~\ref{thm:topology_curves} and~\ref{thm:topology_curves_p1}, wherein we enumerate the path connected components of the moduli space of genus one spectral curves and describe the topology of each component.

\begin{lem}\label{lem:T_graph}
If $p \leq 1$ then there is a diffeomorphism between $\mathcal{\tilde{A}}(p)$, the subset of $\mathcal{\tilde{A}}$ where $S\circ \tilde{\pi} = p$, and
\[
\{(q,k,\tilde{u}) \in \R\times(0,1)\times\R \},
\]
such that fixing a value $q$ gives the level set $\mathcal{\tilde{A}}(p,q)$, defined to be the subset of $\mathcal{\tilde{A}}(p)$ on which we further fix $\tilde{T} = q$. Likewise, if $p \geq 1$ then there is a diffeomorphism between $\tilde{\mathcal{A}}(p)$ and
\[
\{(q,k,\tilde{v}) \in \R\times(0,1)\times\R \},
\]
such that again fixing a value $q$ gives the level set $\mathcal{\tilde{A}}(p,q)$.
In either case, the level sets are graphs over $(0,1)\times \R$.

\begin{proof}
Fix a value of $p$ and consider the function $G(q, k,\tilde{u},\tilde{v}) = \tilde{T}(p,k,\tilde{u},\tilde{v}) - q$ on $\mathcal{\tilde{A}}(p)\times\R$. The preimage $G^{-1}(0)$ is a graph over $\mathcal{\tilde{A}}(p)$ given by $q=\tilde{T}$, so they are diffeomorphic. 
Suppose first that the fixed value of $p\geq 1$. In Lemma~\ref{lem:dTdu_nonzero} we showed that the $\tilde{u}$ derivative of $\tilde{T}$ was non-vanishing and therefore neither is the $\tilde{u}$ derivative of $\tilde{G}$. The implicit function theorem states that there is a function $h$ such that $G^{-1}(0)$ is a graph of the form
\[
\{ (q, k, \tilde{u}, h(q,k,\tilde{u})) \mid q \in \text{Range } \tilde{T}, k \in (0,1), \tilde{u}\in\R \}. 
\]
By Lemma~\ref{lem:range_T}, the range of $\tilde{T}$ is $\R$. This shows that $(q,k,\tilde{u})$ parametrises $\mathcal{A}(p)$ as claimed. If we hold $q$ fixed, then the level set $\mathcal{\tilde{A}}(p,q)$ is parametrised by the remaining two coordinates $(k,\tilde{u})$.

Note the following symmetry:
\begin{align*}
T_0(p,k,u,v) = -p T_0\bra{ \tfrac{1}{p}, k, v, u }.
\labelthis{eqn:T_0 symmetry}
\end{align*}
When $p\leq 1$, this symmetry shows that it is now the $\tilde{v}$-derivative of $\tilde{T}$ that is non-vanishing. Again the implicit function theorem gives the result.
\end{proof}
\end{lem}

With this lemma in hand, we are within striking distance of results about the space of spectral curves $\mathcal{S}_1 \subset \mathcal{A}/\mathbb{Z}_2$. Consider preimage of $\mathcal{S}_1$ in the universal cover $\mathcal{\tilde{A}}$. From Lemma~\ref{lem:closing_conds} the points of the preimage of $\mathcal{S}_1$ are precisely those in which $S\circ\tilde{\pi}$ takes a positive rational value and $\tilde{T}$ takes a positive rational value, we can write this preimage as a disjoint union of level sets
\[
% \mathcal{\tilde{S}} = 
\coprod_{p\in \Q^+,\; q\in \Q} \mathcal{\tilde{A}}(p,q).
% \labelthis{eqn:def tilde S}
\]
We must understand the action of the covering transformations of $\mathcal{\tilde{A}} \to \mathcal{A} \to \mathcal{A}/\mathbb{Z}_2$ as restricted to preimage. The covering transformations of $\mathcal{A}$ over $\mathcal{A}/\mathbb{Z}_2$ are easy to understand; besides the identity there is only $\lambda : (\alpha,\beta) \mapsto (\beta,\alpha)$, which swaps the labelling of the branch points inside the unit disc. 
In terms of the $(p,k,u,v)$ coordinates, $\lambda$ fixes $p = S(\alpha,\beta)$ and $k = k(\alpha,\beta)$, from inspection of equations~\eqref{eqn:def_S} and~\eqref{eqn:def_k}. Recall that $\iu u = f(1)$ and $\iu v = f(-1)$. Our construction relied on a definition of $f$ that sends $\alpha$ to $1$ and $\beta$ to $k^{-1}$. Let $f_s$ be the M\"obius transformation which instead takes $\beta$ to $1$ and $\alpha$ to $k^{-1}$ and consider the composition of $f_s \circ f^{-1}$. It is a M\"obius transformation that exchanges $1$ and $k^{-1}$ and also $-1$ and $-k^{-1}$. It therefore must be the map $z \mapsto (kz)^{-1}$ and hence
% Under this map the point $\iu u$ is taken to $-\iu (ku)^{-1}$ and $\iu v$ is taken to $-\iu (kv)^{-1}$. Thus, under the label-swapping involution $\lambda$,
\[
\lambda: (p,k,u,v) \mapsto \bra{ p,k, -(ku)^{-1}, -(kv)^{-1} }.
\]
To gain a geometric understanding of $\lambda$, we may rescale $u$ in the following manner. Let $U=\sqrt{k} u$ so that $\lambda: U \mapsto -U^{-1}$,
which is a half-rotation of the circle $\RP^1$.

% \maketikzfigure{The circle on the left is $\RP^1$ with four points marked. The map $U \mapsto -U^{-1}$ is applied to produce the circle on the right, with corresponding points coloured the same. We can see that the circle on the right is rotated by half. \label{fig:circle rotation}}{tikz/circle_rotation}

We shall use this rescaled coordinate $U$, and likewise $V = \sqrt{k} v$, to construct coordinates on $\mathcal{\tilde{A}}$ such that the covering transformations are translations. 
% This is motivated by the fact that translating the line $\R$ rotates the quotient $\R/\Z$.
In analogy to Lemma~\ref{lem:mathcal tilde C}, we define coordinates $\tilde{U}$ and $\tilde{V}$ on $\mathcal{\tilde{A}}$ to be $U = \tan \tilde{U}/2$ and $V = \tan \tilde{V}/2$.
% \begin{align*}
% U = \tan \frac{\tilde{U}}{2},       &\quad
%     % U^{-1} = \cot \frac{\tilde{U}}{2},  \\
% V = \tan \frac{\tilde{V}}{2}.       
% % &\quad
%     % V^{-1} = \cot \frac{\tilde{V}}{2}.
% \end{align*}
% Using $U = \sqrt{k}u$, t
The change of coordinates from $\tilde{u}$ to $\tilde{U}$ and from $\tilde{v}$ to $\tilde{V}$ is
\[
\tilde{U} = 2\pi\Wind(\tilde{u}) + 2 \atan \left[ \sqrt{k} \tan \frac{\tilde{u}}{2} \right],
\tilde{V} = 2\pi\Wind(\tilde{v}) + 2 \atan \left[ \sqrt{k} \tan \frac{\tilde{v}}{2} \right].
\]
Observe that these formulae fix multiples of $\Z$.
Hence $\tilde{u}$ and $\tilde{U}$ have the same winding number.
Recall also from Lemma~\ref{lem:mathcal tilde C} that the range of the coordinates $\tilde{u}$ and $\tilde{v}$ is restricted to $\tilde{u} < \tilde{v} < \tilde{u} + 2\pi$.
Since $\tan$ and $\atan$ are increasing functions, it follows that the transformation is order preserving, so $\tilde{U} < \tilde{V} < \tilde{U} + 2\pi$ as well.

% These coordinates $\tilde{U}$ and $\tilde{V}$ can be used interchangeably with $\tilde{u}$ and $\tilde{v}$. Consider for example the coordinates $(p,k,\tilde{u},\tilde{v})$ on $\mathcal{\tilde{A}}$ defined in Lemma~\ref{lem:mathcal tilde C}.
% The determinant of the Jacobian of the change of coordinates to $(p,k,\tilde{U},\tilde{V})$ is
% \[
% \det \mathrm{Jac}\,
% = \begin{vmatrix}
% 1 & 0 & 0 & 0\\
% 0 & 1 & 0 & 0 \\
% 0 & \Partial{\tilde{U}}{k} & \Partial{\tilde{U}}{\tilde{u}} & 0 \\
% 0 & \Partial{\tilde{V}}{k} & 0 & \Partial{\tilde{V}}{\tilde{v}}
% \end{vmatrix}
% = \frac{\sqrt{k}}{\cos^2 \tfrac{\tilde{u}}{2} + k \sin^2 \tfrac{\tilde{u}}{2}} \times \frac{\sqrt{k}}{\cos^2 \tfrac{\tilde{v}}{2} + k \sin^2 \tfrac{\tilde{v}}{2}} \neq 0,
% \]
% so this is indeed a valid change of coordinates. The same is true for the coordinates on $\mathcal{\tilde{A}}(p)$ provided by Lemma~\ref{lem:mathcal tilde C}. For $p \leq 1$, we may change coordinates from $(p,q,k,\tilde{u})$ to $(p,q,k,\tilde{U})$, and for $p \geq 1$ from $(p,q,k,\tilde{v})$ to $(p,q,k,\tilde{V})$.

We can now use the coordinates $(p,k,\tilde{U},\tilde{V})$ on $\mathcal{\tilde{A}}$ to find the group of covering transformations $\mathcal{\tilde{A}} \to \mathcal{A}/\mathbb{Z}_2$.
% We can divide the covering transformations into two types: those that push forward to $\mathcal{A}$ to give the identity and those that push forward to give $\lambda$. 
Recall that the projection of the universal cover $\mathcal{\tilde{A}}$ to $\mathcal{A}$ in part reads $U = \tan \tilde{U}/2$. This is the standard covering of the circle by $\R$, so the subgroup of covering transformations that push forward to the identity are generated by
\[
\tilde{\iota}: (p,k,\tilde{U},\tilde{V}) \mapsto (p,k,\tilde{U} + 2\pi ,\tilde{V} +2\pi). 
\]
Any two covering transformations that push forward to $\lambda$ must differ by a transformation that pushes forward to the identity. If $\tilde{\lambda}$ is any lift of $\lambda$ then the group of covering transformations is generated by $\tilde{\lambda}$ and $\tilde{\iota}$. Take
\[
\tilde{\lambda} : (p,k,\tilde{U},\tilde{V}) \mapsto (p,k, \tilde{U} + \pi, \tilde{V} + \pi).
% \labelthis{eqn:tilde lambda action}
\]
If we apply this transformation twice, we see that $\tilde{\lambda}^2 = \tilde{\iota}$ and so in fact the group of covering transformations $\mathcal{\tilde{A}} \to \mathcal{A}/\mathbb{Z}_2$ is generated by $\tilde{\lambda}$ alone.

If we wish to see the effect of the covering transformations on the preimage of the space of spectral curves, we must determine how the value of $\tilde{T}$ changes when it is precomposed with $\tilde{\lambda}$, since the preimage is a collection of its level sets.

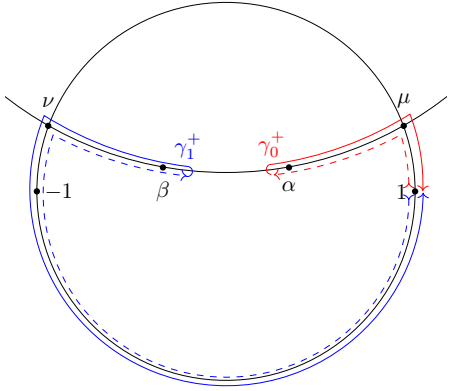
\begin{figure}
\resizebox{0.5\textwidth}{!}{
\begin{tikzpicture}
    \clip (-3.5,-3.5) rectangle (3.5,3.5);
    % \draw[step=0.1,gray,very thin] (-6,-6) grid (6,6);
    
    % unit circle
    \draw (0,0) circle [radius=3];
    % branch point circle
    \draw[color=black] (0,6) circle [radius=5.7];
    % \gamma^+ path, from 1
    \draw[color=red,dashed,>->] (2.9,0.05)
        arc (0:17:2.9)
        arc (-62:-83:5.8);
    \draw[color=red,->] (0.72,0.27)
        arc (-90:-277:0.08)
        node[above, color=red]{$\gamma_0^+$}
        arc (-82:-58:5.6)
        arc (21:-2:3.1);
    
    % \gamma_1^+ path, from 1 but around ν
    \draw[color=blue,dashed,>->] (2.9,-0.05)
        arc (0:-199:2.9)
        arc (-118:-96:5.8);
    \draw[color=blue,->] (-0.64,0.24)
        arc (-110:80:0.08)
        node[above, color=blue]{$\gamma_1^+$}
        arc (263:238:5.6)
        arc (158:359:3.12);
    
    % \gamma^- path, from -1
    % \draw[color=blue,dashed,>->] (-2.9,0)
    %     node[above right, color=blue]{$\gamma^-$}
    %     arc (180:74:2.9)
    %     arc (-25:-35:7.6)
    %     arc (135:315:0.16);
    % \draw[color=blue,->] (0.35,1.40)
    %     arc (325:337:7.9)
    %     arc (68:180:3.1);
    
    % points
    \fill (3,0) circle (0.05) node[left, color=black]{$1$};
    \fill (-3,0) circle (0.05) node[right, color=black]{$-1$};
    \fill (2.82,1.04) circle (0.05) node[above=5pt, color=black]{$\mu$};
    \fill (-2.82,1.04) circle (0.05) node[above=5pt, color=black]{$\nu$};
    \fill (1,0.38) circle (0.05) node[below=3pt, color=black]{$\alpha$};
    \fill (-1,0.38) circle (0.05) node[below=3pt, color=black]{$\beta$};
\end{tikzpicture}
}
\caption{The principal path $\gamma_0^+$ in red and path $\gamma_1^+$ in blue.\label{fig:deck_T_1}}
\end{figure}

\begin{lem}\label{lem:T shift}
The effect of precomposing $\tilde{T}$ with $\tilde{\lambda}$ is to increase its value by $S-1$. That is,
\[
\tilde{T} \circ \tilde{\lambda} - \tilde{T}
= S-1.
\]
\begin{proof}
We first consider the effect of $\lambda$ on $T_0$ in the $\zeta$-plane. Suppose that $\mu$ and $\nu$ are chosen such that $\mu,1,-1$and $\nu$ are arranged clockwise, as shown in Figure~\ref{fig:deck_T_1}. The principal choice of path $\gamma_0^+$, shown in red, traverses around $\alpha$ whereas $\gamma_1^+$ is the principal choice after swapping the labels of the roots, shown in blue.
The difference between these two paths is homologous to $-B$. 
% So by the construction of $\Theta^P$,
% \[
% \int_{\gamma_1^+} \Theta^P - \int_{\gamma_0^+} \Theta^P = \int_{\S^1} \Theta^P = -2\pi\iu.
% \]
The same is true for $\gamma_0^-$.
% Likewise if we consider the difference between the principal path $\gamma_0^-$ and the path $\gamma_1^- $,  we again have a anticlockwise loop of the upper unit circle.
% , shown in Figure~\ref{fig:deck_T_-1}.
% Hence
% \[
% \int_{\gamma_1^-} \Theta^P - \int_{\gamma_0^-} \Theta^P = \int_{\S^1} \Theta^P = -2\pi\iu.
% \]
Putting these together we conclude that the value of $T_0$ changes by $1-S$ under this transformation at points where $\mu,1,-1$and $\nu$ are so ordered:
\begin{align*}
2\pi\iu (T_0 \circ \lambda - T_0)
&= \bra{ S\int_{\gamma_1^-} \Theta^P - \int_{\gamma_1^+} \Theta^P } - \bra{ S\int_{\gamma_0^-} \Theta^P - \int_{\gamma_0^+} \Theta^P } \\
&= S(-2\pi\iu) - (-2\pi\iu) \\
&= 2\pi\iu (1-S).
\end{align*}
To infer the effect of the transformation on $\tilde{T}$ however we must take into account how the coordinates $\tilde{u}$ and $\tilde{v}$ may have changed their winding numbers. 
As one traverses the unit circle clockwise in the $\zeta$-plane, one traverses the imaginary axis in the $z$-plane upwards. 
Hence $0 = f(\mu) < u < v < \infty = f(\nu)$.
% then this restricts the arrangement of $\iu u$ and $\iu v$. By definition,
% \[
% 0 = f(\mu),\;\; \iu u = f(1),\;\; \iu v = f(-1),\;\; \infty = f(\nu),
% \]

As $u$ and $v$ are both positive, it must be that $\tilde{u},\tilde{v} \in (2\pi n, 2\pi n + \pi)$. Hence $\tilde{U}$ and $\tilde{V}$ also lie in that interval. Under applying the transformation $\tilde{\lambda}$, the coordinates $\tilde{U}$ and $\tilde{V}$ will be translated by $\pi$, increasing their winding numbers by $1$. Combining the effect of $\lambda$ on $T_0$ with the change of winding number in~\eqref{eqn:tilde T computable} shows the relation on this particular subset of $\mathcal{\tilde{A}}$. But a relation between analytic functions on an open set applies everywhere.
\end{proof}
\end{lem}

Thus we can describe the effect of the covering transformation $\tilde{\lambda}$ on points in the level set $\mathcal{\tilde{A}}(p,q)$, where $\tilde{T} = q$. Take a point $(q,k,\tilde{X}) \in \mathcal{\tilde{A}}(p)$, where if $p\leq 1$ let $\tilde{X}$ be $\tilde{U}$ and otherwise take it to be $\tilde{V}$ (cf. Lemma~\ref{lem:T_graph}). Under the covering transformation $\tilde{\lambda}$,
\[
\bra{p,q,k,\tilde{X}} \mapsto \bra{p, q + (p-1), k, \tilde{X} + \pi}.
\labelthis{eqn:group action}
\]
Viewing the cosets of the group of covering transformations as an equivalence relation, for $p\neq 1$ they provide an identification between the level sets $\mathcal{\tilde{A}}(p,q)$ and $\mathcal{\tilde{A}}(p,q + l(p-1))$ for any integer $l$.

\begin{thm}\label{thm:topology_curves}
For $p\neq 1$, the subspace $\mathcal{A}(p)/\mathbb{Z}_2$ of $\mathcal{A}/\mathbb{Z}_2$, where $S=p$, is diffeomorphic to
\[
\left\{ \bra{[q],k,\tilde{X}} \in \bra{\R/ (p-1)\Z} \times (0,1) \times \R \right\}.
\]
If $p\in\Q$ then the subspace of spectral curves is obtained by restricting $q$ to be rational;
\[
\mathcal{S}_1 \cap \mathcal{A}(p)/\mathbb{Z}_2  = \left\{ \bra{[q],k,\tilde{X}} \mid [q] \in \Q/ (p-1)\Z \right\}.
\]
\begin{proof}
Fix $p\neq 1$ and consider $\mathcal{A}(p)/\mathbb{Z}_2$. It is the quotient of $\mathcal{\tilde{A}}(p)$ by the group of covering transformations $\Z\langle\tilde{\lambda}\rangle$. By Lemma~\ref{lem:T_graph}, $\mathcal{\tilde{A}}(p)$ is foliated by $\mathcal{\tilde{A}}(p,q)$ and we have just shown in~\eqref{eqn:group action} how different $\mathcal{\tilde{A}}(p,q)$ can be identified if their values of $q$ differ by a multiple of $p-1$. Hence it is sufficient to take one representative from each element of $\R/(p-1)\Z$ to cover the image.
Further, from Lemma~\ref{lem:closing_conds} we know that a point of $\mathcal{A}/\mathbb{Z}_2$ corresponds to a spectral curve exactly when $S=p$ and $T=q$ are rationally valued.
\end{proof}
\end{thm}

This leaves just one special case, where $p=1$. We see that the action of $\tilde{\lambda}$ on $\mathcal{\tilde{A}}(1)$ fixes the value of $\tilde{T}$ and so does not identify different level sets $q=\tilde{T}$. Instead, the group action creates an equivalence relation on each level set. Specifically, the action of $\tilde{\lambda}$ on $\mathcal{\tilde{A}}(1,q)$ given by~\eqref{eqn:group action} reads $\bra{1,q,k,\tilde{X}} \mapsto \bra{1, q, k, \tilde{X} + \pi}$.

\begin{thm}\label{thm:topology_curves_p1}
The space $\mathcal{A}(1)/\mathbb{Z}_2$ is diffeomorphic to
\[
\left\{ \bra{q,k,\left[\tilde{X}\right]} \in \R \times (0,1) \times \R/\pi\Z \right\},
\]
the product of $\R$ and an annulus, such that the subspace of spectral curves is the restriction of the first component of the product to $\Q$.

\begin{proof}
The orbit of a point $(1,q,k,\tilde{X})$ of $\mathcal{\tilde{A}}(1)$ is $\bra{1,q,k,\tilde{X}+\pi\Z}$,
so the quotient sends $\tilde{X} \in \R$ to $\left[\tilde{X}\right] \in \R/\pi\Z$. The conclusion about the spectral curves again follows immediately from Lemma~\ref{lem:closing_conds}.
\end{proof}
\end{thm}

These two theorems deliver an understanding of the subspace of spectral curves $\mathcal{S}_1$ within $\mathcal{A}/\mathbb{Z}_2$. 
For $p$ not equal to one, the components of $\mathcal{S}_1(p)$ are contractible and dense within $\mathcal{A}(p)/\mathbb{Z}_2$. As can be seen in Figure~\ref{fig:p05 plot} below, they should be thought of as being intertwined around a central axis $\Delta\cap\{S=p\}$, similar to a family of helicoids. But for $p$ equal to one, instead we have a dense collection of annuli $\mathcal{S}_1(1)$, as in Figure~\ref{fig:p1 plot}.

Unlike true helicoids however, cross-sections perpendicular to the central axis do not meet the axis. Instead they spiral infinitely closer. In this respect, they behave like the cone of a spiral. This aspect is especially prominent in Figure~\ref{fig:p1 plot}, which one could think of as a family of cones with a common vertex (though the plotting software has difficulty for $k\approx 1$, so the `cones' are truncated and appear to be heading towards a white hole in the middle of the figure). Some aspects of this vertex structure are known to the authors and will be addressed in a future paper.
We saw in the introduction to this chapter that the parameter $p$ can be thought of as controlling the slope of the level sets as they wind around the central axis. From this point of view, $p=1$ is the intermediate case between right and left handed spirals where the slope is `flat' and the level sets `close up'.

\begin{figure}[p]
    \centering
    \includegraphics[height=0.4\textheight]{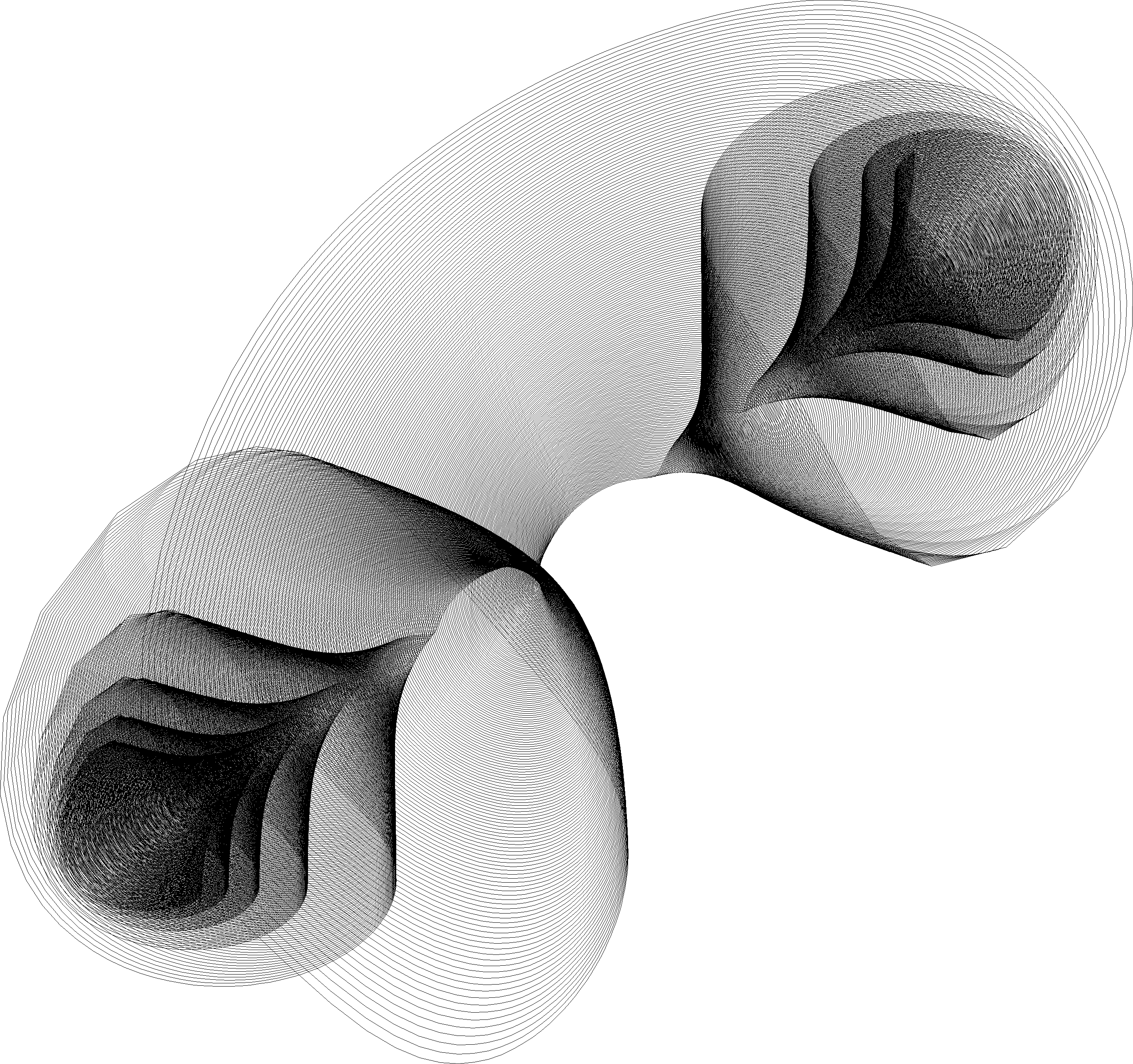}
    \caption{One component of the moduli space $\mathcal{S}_1(0.5)$. Notice that the `helicoid' is wrapped in a left handed direction around a central axis.\label{fig:p05 plot}}

    \vspace*{\floatsep}

    \includegraphics[height=0.37\textheight]{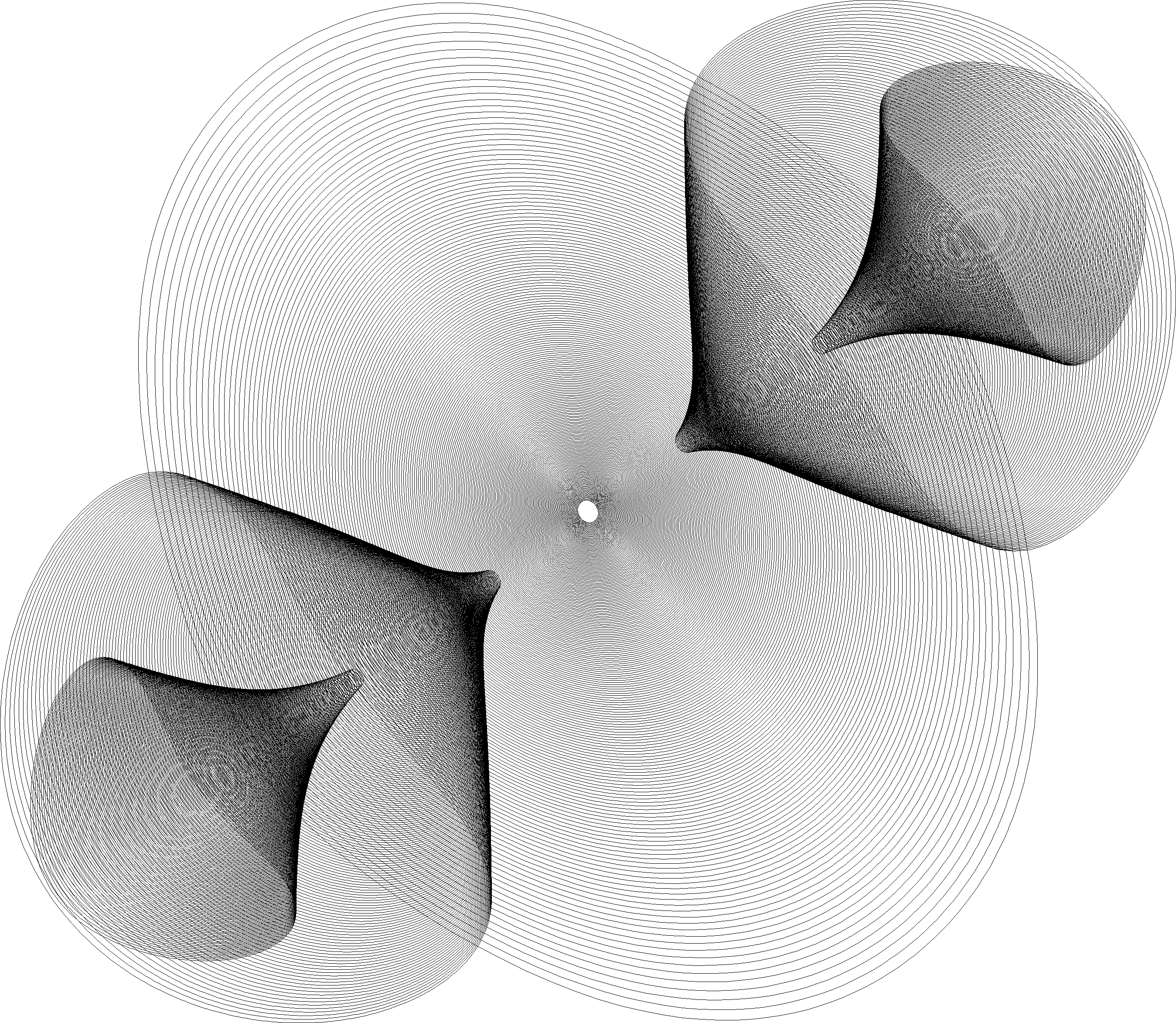}
    \caption{Select components of the moduli space $\mathcal{S}_1(1)$, namely those components on which $T$ is $-3$, $-1$, $0$, $1$, or $3$. The component on which $T$ is $0$ is the disc in the middle. If one is viewing this article digitally, one may like to zoom these images.\label{fig:p1 plot}}
\end{figure}

% \begin{figure}
%     \includegraphics[width=\textwidth]{graphics/moduli_plot_p2.png}
%     \caption{One component of the moduli space $\mathcal{S}_1(2)$. Notice, contra Figure~\ref{fig:p05 plot}, that the `helicoid' is wrapped in a right handed direction around a central axis.\label{fig:p2 plot}}
% \end{figure}

\subsection{The Moduli Space of Spectral Data}
In this subsection we use our knowledge of the space of spectral curves $\mathcal{S}_1$ to describe the space of spectral data $(\Sigma,\Theta^1,\Theta^2,E)\in\mathcal{M}_1$. We treat $\mathcal{M}_1$ as a bundle over $\mathcal{S}_1$ via projection to the first element of the tuple.

First let us establish the local structure of $\mathcal{M}_1$. 
Fix a path connected component $\mathcal{X}$ of $\mathcal{S}_1$.
The function $S$ is well defined on $\mathcal{X}$ but $T$ may not be. 
However, around any point take a simply connected neighbourhood $\mathcal{V}$. On $\mathcal{V}$ it is possible to choose paths $\gamma^+$ and $\gamma^-$ on each spectral curve that vary smoothly with changes of the branch points. Hence the choice of paths $\gamma^+$ and $\gamma^-$ can be made continuously on $\mathcal{V}$ and Lemma~\ref{lem:closing_conds} gives us well-defined and smoothly varying differentials $\Psi^E$ and $\Psi^P$. The space of differentials satisfying~\ref{P:poles}--\ref{P:closing} is then trivial $\Z^2$-bundle over $\mathcal{V}$ with basis $\Psi^E$ and $\Psi^P$.

We recall from Section~\ref{sec:Genus Zero} the definition of the integer matrices 
\[
\Mat_2^*\Z = \Set{ M \in \Mat_2\Z }{\det M \neq 0}. 
\]
There the moduli space $\mathcal{M}_0$ of spectral data with a genus zero spectral curve was described as the product $\mathcal{M}_0 = D \times \Mat_2^*\Z$. Similarly, we may describe the space of triples $(\Sigma,\Theta^1,\Theta^2)$ of spectral data whose spectral curve lies in $\mathcal{V}$ as $\mathcal{V} \times \Mat_2^*\Z$ where the differentials $(\Theta^1,\Theta^2)$ are given by the matrix and the frame $(\Psi^E,\Psi^P)$ as
\[
\begin{pmatrix}
\Theta^1 \\ \Theta^2
\end{pmatrix}
=
\begin{pmatrix}
b_1 & l_1 \\
b_2 & l_2
\end{pmatrix}
\begin{pmatrix}
\Psi^E \\ \Psi^P
\end{pmatrix}.
\]
The non-vanishing of the determinant of this matrix is equivalent to Condition~\ref{P:linear independence}, that the pair of differentials $\Theta^1, \Theta^2$ are linearly independent. As discussed in~\cite[p.~690]{Hitchin1990}, the line bundle $E$ may be chosen in this case from a one dimensional subvariety of the Jacobian; topologically this is a circle.

With the local structure established, it remains to show the global structure of $\mathcal{M}_1$. From the previous section, the path connected components of $\mathcal{S}_1$ come in two types: either $\mathcal{S}_1(p,[q])$ for $p\in\Q^+, p\neq 1$ and $[q]\in\Q/p\Z$, or  $\mathcal{S}_1(1,q)$ for $q\in\Q$. We similarly decompose $\mathcal{M}_1$ according to which component its spectral curve belongs. As the components $\mathcal{S}_1(p,[q])$ are contractible, we have immediately that
\[
\mathcal{M}_1(p, [q]) = \mathcal{S}_1(p,[q]) \times \Mat_2^*\Z \times \S^1.
\]

For $p=1$, let us fix a path component $\mathcal{S}_1(1,q)$.
This component is an annulus and not simply connected. The differential $\Psi^P$ is not well-defined on all of $\mathcal{S}_1(1,q)$ because it is not possible to make a consistent continuous choice of paths $\gamma^+$ and $\gamma^-$. For example, consider the family of spectral curves in $\mathcal{S}_1(1,0)$ with branch points at $0.5 e^{\iu t}$ and $-0.5 e^{\iu t}$. If one takes $\gamma^+$ to be a path that wraps around $\zeta = 0.5$ on $\Sigma_{t=0}$ then as $t$ is increased to $\pi$, the same path now passes around the other branch point $\zeta=-0.5$ on the same spectral curve $\Sigma_{t=\pi} = \Sigma_{t=0}$. 

We must therefore describe the monodromy of $\Psi^P$ along a loop $\ell : [0,1] \to \mathcal{S}_1(1,q)$ that wraps around the annulus once. 
Choose $\gamma^+(t)$ and $\gamma^-(t)$ on each $\Sigma_{\ell(t)}$ such that the paths vary smoothly in $t$.
%  then note in particular that $\gamma^+(0)$ and $\gamma^+(1)$ will be different paths on the same curve. One may ask: if we construct $\Psi^P$ using $\gamma^+(t)$ and $\gamma^-(t)$ then what will be the change in $\Psi^P$ when we return to $\ell(1) = \ell(0)$?
From~\eqref{eqn:group action} it follows that $T = n'/m'$ is well defined and constant on the whole annulus. We may simplify so that it is a reduced fraction. Writing $S=1=n/m$ as a reduced fraction implies that $m=n=1$. The period of $\Psi^P$ is therefore $l = m' / \gcd(m',mn') = m'$ and we may write $\Psi^P(0) = b(0)\Theta^E + m' \Theta^P$ for the differential on $\Sigma_{\ell(0)}$. As the periods of the differential are integral, they cannot change along the path $\ell$ and so $\Psi^P(t) = b(t) \Theta^E + m' \Theta^P$.
As we move from the start of $\ell$ to the end, from our previous computation in Lemma~\ref{lem:T shift}, each of
\[
\int_{\gamma^+(t)} \Theta^P \;\text{ and } \int_{\gamma^-(t)} \Theta^P
\]
will be incremented or decremented (depending on the orientation of $\ell$) by $2\pi\iu$. The difference $\Psi^P(1) - \Psi^P(0) = (b(1)-b(0))\Theta^E$ can be explicitly computed from~\eqref{eqn:def Psi coeff} as
\begin{align*}
b(1)
&= \frac{1}{2\iu \eta^+(1)}\bra{ 2\pi\iu y n'm - m' \bra{\int_{\gamma^+(1)} \Theta^P} }
% &= \frac{1}{2\iu \eta^+(1)}\bra{ 2\pi\iu y n'm - m' \bra{\int_{\gamma^+(0)} \Theta^P + 2\pi\iu} } \\
% &= b(0) - \frac{2\pi\iu}{2\iu \eta^+(1)}m' \\
= b(0) - am',
\end{align*}
so that $\Psi^P(1) - \Psi^P(0) = - m' \Psi^E$. This shows that every loop around the annulus shifts $\Psi^P$ by $m' \Psi^E$.

This causes the space of triples $(\Sigma,\Theta^1,\Theta^2)$ such that $\Sigma \in \mathcal{S}_1(1,q)$ to be simply connected. Given any closed path $\ell$ in it, we may project this loop down to $\mathcal{S}_1(1,q)$. 
Suppose for contradiction that the projection is non-trivial, i.e.\ that it winds a nonzero number of times around the annulus. By assumption the differentials are unchanged from the beginning to end of the path, so either $m' = 0$ or $l_1=l_2 = 0$. The former is excluded because by definition $n'/m' = q$ and so $m'$ is never zero. The latter implies that the differentials are both multiples of $\Theta^E$, which by~\ref{P:linear independence} contradicts their linear independence. 
Of course, if the projection of the loop is null-homotopic, then it is contained in a simply connected neighbourhood $\mathcal{V}$ of $\mathcal{S}_1(1,q)$ and we may use the frame $\langle \Psi^E,\Psi^P \rangle$ to lift to a null-homotopy of $\ell$ in the space of triples.
Hence each connected component is diffeomorphic to the universal cover of the annulus $\mathcal{S}_1(1,q)$, a ribbon $(0,1)\times \R$. This is contractible, so it follows that $\mathcal{M}_1$ is trivial as am $\S^1$-bundle over the space of triples. 

We can translate the shift of $\Psi^P$ into a statement about $\Mat_2^*\Z$. Given a tuple of spectral data $(\Sigma,\Theta^1,\Theta^2)$ in $\mathcal{M}_1(1,q)$, one is free to vary $\Sigma$ within $\mathcal{S}_1(1,q)$, but the effect on the differentials of making loop around  the annulus is
\begin{align*}
\begin{pmatrix}
\Theta^1 \\ \Theta^2
\end{pmatrix}
% &=
% \begin{pmatrix}
% b_1 & l_1 \\
% b_2 & l_2
% \end{pmatrix}
% \begin{pmatrix}
% \Psi^E \\ \Psi^P
% \end{pmatrix} \\
\mapsto
\begin{pmatrix}
b_1 & l_1 \\
b_2 & l_2
\end{pmatrix}
\begin{pmatrix}
    \Psi^E \\ \Psi^P - m' \Psi^E
\end{pmatrix}
=
\begin{pmatrix}
b_1 & l_1 \\
b_2 & l_2
\end{pmatrix}
\begin{pmatrix}
    1 & 0 \\
    -m' & 1
\end{pmatrix}
\begin{pmatrix}
    \Psi^E \\ \Psi^P
    \end{pmatrix}.
\end{align*}
If $B_q$ is the subgroup of $\Mat_2^*\Z$ of matrices of the form
\[
\begin{pmatrix}
1 & 0 \\
m'\Z & 1
\end{pmatrix},
\]
then the connected components of $\mathcal{M}_1(1,q)$ are enumerated by the right $B_q$-orbits of $\Mat_2^*\Z$.

In summary, $\mathcal{M}_1$ is the disjoint union
\begin{gather*}
\mathcal{M}_1
% =
% \coprod_{q \in \Q} \mathcal{M}_1(1,q)
% \;\;\amalg \coprod_{\substack{p \in \Q^+,\, p \neq 1\\ [q] \in \Q/(p-1)\Z}}
% \mathcal{M}_1(p,[q]) \\
=
\coprod_{q \in \Q} \Big[ \widetilde{\mathcal{S}_1(1,q)} \times \bra{\Mat_2^*\Z / B_q} \times \S^1 \Big]
\;\amalg
\coprod_{\substack{p \in \Q^+,\, p \neq 1\\ [q] \in \Q/(p-1)\Z}} \Big[ \mathcal{S}_1(p,[q]) \times \Mat_2^*\Z \times \S^1 \Big] \\
=
(0,1)\times\R \times \left[
\coprod_{q \in \Q} \Mat_2^*\Z / B_q
\;\amalg
\coprod_{\substack{p \in \Q^+,\, p \neq 1\\ [q] \in \Q/(p-1)\Z}} \Mat_2^*\Z
\right] \times \S^1.
\end{gather*}

\section{Corollaries}\label{sec:Corollaries}

In this final section, we shall extend the symmetry exhibited in~\eqref{eqn:T_0 symmetry} to a general transformation $\chi$ on the space of spectral curves of a fixed genus $\mathcal{S}_g$ and give a geometric interpretation. We will illustrate this interpretation by observing how it applies to the space $\mathcal{S}_0$ of spectral curves with genus zero. 
Then we will examine the special case of harmonic maps to a $2$-sphere and show that those with a genus one spectral curve can be identified with a particular path component of $\mathcal{S}_1(1)$.

Recall the symmetry in~\eqref{eqn:T_0 symmetry},
\[
T_0(p,k,u,v) = -p T_0\bra{ \tfrac{1}{p}, k, v, u }.
\]
How should one interpret it geometrically? Looking at the transformation
\[
p = \frac{\abs{1-\alpha}\abs{1-\beta}}{\abs{1+\alpha}\abs{1+\beta}}
\mapsto \frac{1}{p} = \frac{\abs{1+\alpha}\abs{1+\beta}}{\abs{1-\alpha}\abs{1-\beta}},
\]
the natural hypothesis is that it is induced by $\chi:(\alpha,\beta)\mapsto (-\alpha,-\beta)$.
% \begin{align*}
% \chi: \mathcal{A} &\to \mathcal{A} \\
% (\alpha,\beta) &\mapsto (-\alpha,-\beta).
% % \labelthis{eqn:def chi}
% \end{align*}
Indeed this can be seen to be the case, as $k$ is invariant under such a transformation and the associated map between the hyperelliptic curves
\[
\chi_{(\alpha,\beta)}: \Sigma(\alpha,\beta) \to \Sigma(-\alpha,-\beta),
\quad (\zeta, \eta) \mapsto (-\zeta,-\eta)
\]
interchanges $1$ and $-1$, in effect swapping the roles of $u = -\iu f(1)$ and $v = -\iu f(-1)$. The pullback of the differentials under the map $\chi_{(\alpha,\beta)}$ preserves the integrality of the periods and the closing conditions. Because it commutes with $\rho$ and $\sigma$, the symmetry properties of the differentials and line bundle are also preserved. Hence spectral data on $\Sigma(\alpha,\beta)$ is transformed into spectral data on $\Sigma(-\alpha,-\beta)$.
The same reasoning applies in general to spectral curves of any genus.
The harmonic map $g(z) : \mathbb{T}^2 \to \S^3$ arises from the spectral data as the gauge transformation between the connections corresponding to $\zeta=1$ and $\zeta=-1$ in~\eqref{eqn:flat connection translation}, so exchanging these points with $\chi_{(\alpha_1,\alpha_2,\dots,\alpha_g)}$ gives the inverted map $g(z)^{-1}$, which is also harmonic~\cite[Prop~8.2]{Uhlenbeck1989}.

This can be seen explicitly in the genus zero case, which was treated in Section~\ref{sec:Genus Zero}. Recall that the equation of any harmonic map corresponding to spectral data with a genus zero spectral curve may, as in~\eqref{eqn:genus zero simple map}, be written
\[
g(w_R + \iu w_I) = \exp (-4w_R X) \exp (4w_I Y),
\]
for
\[
X = \norm{X}\begin{pmatrix}
0 & 1 \\
-1 & 0
\end{pmatrix}, \quad
Y = \norm{Y}\begin{pmatrix}
0 & e^{\iu \delta} \\
-e^{-\iu \delta} & 0
\end{pmatrix}.
\]
The norms, $\norm{X}$ and $\norm{Y}$, come from the identification of $\su_2$ with $\R^3$. Geometrically, $\delta$ is the angle between $X$ and $Y$. Under inversion,
\[
g(w)^{-1} = \exp (-4w_I Y) \exp (4w_R X).
\]
To bring this back into the form of~\eqref{eqn:genus zero simple map}, we must perform two operations. First, we must change coordinates on the domain so that the real part is in the first factor. 
Second we must rotate the image so that the matrices $-Y$ and $X$ take the required forms. After this restoration, we find that the angle parameter has become $\pi-\delta$. As an aside, 
% as seen in Figure~\ref{fig:genus0 linked} 
this angle parameter determines the image of the harmonic map up to an $SO(4)$ rotation of $\S^3$. In particular, these maps $g$ and $g^{-1}$ have congruent images.

Recall that the parameter $x$ of $g$ was defined by~\eqref{eqn:conformal type} to be $\norm{Y}/\norm{X}$. The value of this ratio is inverted for $g^{-1}$. Now using~\eqref{eqn:def branch point genus zero} to determine the branch point of the spectral curve associated to $g^{-1}$ we have
\[
\frac{\frac{1}{x} e^{\iu (\pi-\delta)} - \iu}{\frac{1}{x} e^{\iu (\pi-\delta)} + \iu}
% = \frac{-e^{-\iu \delta} - \iu x}{-e^{-\iu \delta} + \iu x}
% = \frac{-1 - \iu x e^{\iu \delta}}{-1 + \iu x e^{\iu \delta}}
= -\frac{x e^{\iu \delta}-\iu}{x e^{\iu \delta} + \iu} 
= -\alpha,
\]
which is to say the branch point of the spectral curve associated to $g^{-1}$ is the negative of the branch point of the spectral curve associated to $g$, as asserted above.

Returning to genus one spectral curves, of special interest are the spectral curves that are a fixed point of this transformation $\chi: \Sigma(\alpha,\beta) \mapsto \Sigma(-\alpha,-\beta)$. For these spectral curves $\chi_{(\alpha,\beta)}$ is an extra involution and $\beta=-\alpha$. If we can show that they admit spectral data these would be exactly the genus one curves that would meet the conditions for the associated harmonic map to have a totally geodesic two-sphere as its image~\cite[Theorem~8.20]{Hitchin1990}. Hitchin~\cite[p693]{Hitchin1990} identifies a particular one parameter family of these maps as the Gauss maps of Delaunay surfaces. Burstall and Kilian~\cite{Burstall2006} also examine these Delaunay surfaces and their associated families from a loop-group perspective. We shall show that all such curves in fact admit spectral data and further identify in which component of $\mathcal{S}_1$ they reside.

Suppose that $\Sigma(\alpha,-\alpha)$ is a curve that is fixed by $\chi$. As $\chi$ sends $p \mapsto p^{-1}$, it follows that $p$ is $1$. To show that $\Sigma(\alpha,-\alpha)$ admits spectral data it remains to shows that $T$ is rationally valued at this point $(\alpha,-\alpha) \in \mathcal{A}$. But note that $\chi$ and $\lambda$ act in the same way on $\{(\alpha,-\alpha) \in \mathcal{A}$, as one would expect because they are both swapping $\alpha$ and $-\alpha$, the two branch points inside the unit disc. Precomposing $T_0$ with $\lambda$ shifts its value by $1-p$, which in this case fixes $T_0$. On the other hand precomposing $T_0$ with $\chi$ negates the function. Together this shows that $T_0$ is zero for these curves. Conversely, the disjoint annuli that constitute $\mathcal{S}_1(1)$ are determined uniquely by the value $T_0$ takes on them.
Therefore we have shown that $\{\Sigma(\alpha,-\alpha) \in \mathcal{A}/\mathbb{Z}_2\}$ is the only annulus of $\mathcal{S}_1(1)$ where $T_0$ is zero, and these are exactly the harmonic maps to the sphere with a genus one spectral curve.

\bibliographystyle{ross_bibsty}
\bibliography{zotero_real}

\vfill
~\\
Dr Emma Carberry, \\
School of Mathematics and Statistics, University of Sydney, NSW 2006, Australia \\
email: \texttt{emma.carberry@sydney.edu.au}\\
~\\
Dr Ross Ogilvie, \\
Mathematics Chair III, Universit\"at Mannheim, D-68131 Mannheim, Germany \\
email: \texttt{r.ogilvie@math.uni-mannheim.de}

\end{document}